\documentclass[11pt]{article}
\usepackage{amsfonts,latexsym,graphics}
\newtheorem{lemma}{Lemma}
\newtheorem{theorem}{Theorem}
\newtheorem{definition}{Definition}
\newtheorem{corollary}{Corollary}
\newtheorem{proposition}{Proposition}
\def\bideg{\mathrm{bideg}\,}
\def\tdeg{\mathrm{tdeg}\,}
\def\sdeg{\mathrm{sdeg}\,}

\def\mod{\mathrm{mod}\,}
\def\In{\mathrm{In}\,}
\def\Im{\mathrm{Im}\,}
\def\Ker{\mathrm{Ker}\,}
\def\Hom{\mathrm{Hom}\,}
\def\Tor{\mathrm{Tor}\,}
\def\rk{\mathrm{rk}\,}
\def\codim{\mathrm{codim}\,}

\def\R{{\mathbb R}}
\def\C{{\mathbb C}}
\def\Z{{\mathbb Z}}
\def\A{{\mathcal A}}
\def\B{{\mathcal B}}
\def\K{{\mathcal K}}
\def\M{{\mathcal M}}

\def\E{{\mathcal E}}
\def\G{{\mathcal G}}
\def\D{{\mathcal D}}
\def\F{{\mathcal F}}
\def\L{{\mathcal L}}
\def\V{{\mathcal V}}
\def\H{{\mathcal H}}

\begin{document}

\vskip12mm
\begin{center}
{\Large Massey products in symplectic manifolds}
\footnote{Supported by Russian Foundation for Basic Researches
(the grants 99-01-01202-a (I. K. B.),
96-15-96877 and 98-01-00749 (I. A. T.)}
\end{center}

\begin{center}
{\large Ivan K. BABENKO and Iskander A. TAIMANOV}
\end{center}

\date{}

\vskip7mm

\begin{center}
{\bf \S 1. Introduction}
\end{center}

In this paper we resume our study of the formality problem for 
symplectic manifolds, which we started in \cite{BT0,BT} 
where the first examples of nonformal simply connected 
symplectic compact manifolds were constructed and a method of constructing 
such manifolds by symplectic blow-ups was introduced.

A smooth manifold $X$ is called symplectic if there is a nondegenerate closed
$2$-form $\omega$ on $X$. We denote this object by $(X,\omega)$. 
Symplectic manifolds first appeared in analytic mechanics but
long ago became a subject of pure mathematics \cite{McDS}. 
In this paper we consider only closed symplectic manifolds.

It is naturally to expect that the existence of a symplectic structure 
strongly restricts the topology of a manifold. A lack of a wide spectrum
of examples only confirms this guess. The main examples of symplectic 
manifolds are K\"ahler manifolds and the existence of a K\"ahler structure 
really poses strong topological restrictions on a manifold.

The problem of existence of symplectic manifolds with no K\"ahler structure
was posed in the early 70s and the first example of such a manifold was
a $4$-dimensional non-simply connected symplectic manifold 
constructed by Thurston \cite{Th}. It appeared
later that this manifold is one of complex non-K\"ahler surfaces from the
Kodaira classification list. Using this example McDuff constructed a simply
connected symplectic manifold with no K\"ahler structure \cite{McD}.
This implies that the class of K\"ahler manifolds is a subclass of
symplectic manifold but the problem how narrow is it is far from a 
complete solution.

Recently Gompf showed that any finitely presented group is the fundamental 
group of some $4$-dimensional symplectic manifold \cite{G}.
Notice that the first Betti number of a K\"ahler manifold is even.

One of the basic topological properties of K\"ahler manifolds is a formality
established in \cite{DGMS}. The formality of a space means that its real
homotopy type is completely defined by the real cohomology ring
$H^{\ast}(X,\R)$.

An existence of nonformal and, hence, non-K\"ahler symplectic manifolds 
was 
known quite long ago \cite{CFG}. In particular, the Kodaira--Thurston
manifold is nonformal although Thurston did not discuss that.
All examples of nonformal symplectic manifolds known until recently were
non-simply connected and their nonformality arises from certain properties
of the fundamental group.  For simply connected manifolds 
there was the conjecture by Lupton and Oprea, which stated that simply
connected symplectic manifolds are formal and in this sense are similar to 
K\"ahler manifolds \cite{LN1661}. Some partial results confirming that
were obtained in \cite{LO}.

This conjecture was disproved by us in \cite{BT0,BT} where infinitely
many nonformal simply connected symplectic manifolds of all even 
dimensions 
greater than $8$ were presented.
In particular, the $10$-dimensional McDuff manifold \cite{McD} is nonformal.
Any simply connected manifold of dimension less than equal $6$ is formal
\cite{Mil,NM}. Therefore the problem on existence of nonformal simply 
connected symplectic manifolds is open only for $8$-dimensional manifolds.

A well-known obstruction to formality of a space $X$ is an existence of 
nontrivial higher, i.e., $n$-tuple for $n \geq 3$, Massey products in 
$H^{\ast}(X)$. It is a nontriviality of certain triple ordinary
Massey products in some simply connected symplectic manifolds that
was established in \cite{BT}.

Independently of the formality problem an existence of higher Massey products 
reflects a complexity of the topology of a space. Higher Massey products, or 
operations, appear in some problems of homotopy topology, which
are, for instance, a description of the cohomology Hurewicz homomorphism 
and a description of the differentials of the Eilenberg--Moore spectral 
sequence. One of the main aims of this paper is to show that symplectic and,
in particular, simply connected symplectic manifolds may have nontrivial 
$n$-tuple Massey products for arbitrary large $n$.

The paper comprises two parts.

In the first part we introduce the most general Massey products.
A particular case of them are the classical matrix 
Massey products introduced by May in 1969. Our approach to the 
definition of Massey products is based on
considering solutions to the Maurer--Cartan equation 
\begin{equation}
d\,A-\bar{A}\cdot A\equiv 0\ \ \mod K(A).
\label{1}
\end{equation}
on the algebra of
matrix differential forms on the manifold. 
Here $A$ is an upper triangular matrix whose entries are differential forms on
the manifold and $K(A)$ is a certain submodule associated with this 
matrix.

It was noticed by May that the defining system for any matrix Massey product
can be presented by a matrix satisfying the Maurer--Cartan equation
(\ref{1}) but this approach was not developed.

A solution $A$ to the equation (\ref{1}) may be treated as a formal 
connection. Then the corresponding generalized Massey products are the 
cohomology classes of ``the curvature form''
$$
\mu(A)=d\,A-\bar{A}\cdot A
$$
of the connection $A$. This differential-geometric approach to the 
definition of Massey products enables us to establish 
the basic properties of these products in a quite simple manner.

In the first part we consider in brief relations of Massey products to
the Hurewicz and the suspension homomorphisms and to the 
Eilenberg--Moore spectral sequence.

In the second part of this paper we consider applications of 
Massey products to symplectic manifolds. There are three main methods of
constructing new symplectic manifolds from given ones. These are
the symplectic fibration, the fiber connected sum, and the 
symplectic blow-up. Although for constructing symplectic manifolds with
nontrivial Massey products one may use the symplectic fibration 
as it was done, for instance, in the case of the Kodaira--Thurston manifold,
it seems that for these reasons the symplectic blow-up \cite{Gromov1,McD}
is the most perspective construction.

Let $(X,\omega)$ be a symplectic manifold and $Y \subset X$ be its
symplectic submanifold. We denote by $\widetilde{X}$ the symplectic blow-up
of $X$ along $Y$. The main problem studied in the second part is what may 
happen with higher Massey products under the symplectic blow-up.
It splits into two parts:

1) "{\it A survival}". Let $X$ have irreducible higher (matrix) Massey
products. What conditions on $Y$ guarantee that these products survive in
$\widetilde{X}$?

2) "{\it An inheritance}". Assume that $Y$ has nontrivial higher
(matrix) Massey products. Under what conditions on $X$ these products are
inherited by $\widetilde{X}$?

Let us formulate the main results in this direction.

\smallskip

{\bf Theorem A} {\it Let a simply connected symplectic manifold $X$ have 
an irreducible generalized Massey product of dimension $k$. Then for any
symplectic submanifold $Y \subset X$ with $\codim Y>k$ 
the corresponding symplectic 
blow-up $\widetilde{X}$ also has an irreducible generalized Massey product of
dimension $k$.}

\smallskip

{\bf Theorem B} {\it Let a symplectic manifold $(Y,\omega)$ have a nontrivial 
matrix $n$-tuple Massey product $\langle S_1,\dots,S_n\rangle$ where
$S_i\in N(H^1(Y))$ are matrices of one-dimensional
cohomology classes for $1\leq i\leq n$. 
Then for any symplectic embedding $Y\subset X$ of codimension
not less than $2(n+1)$ the corresponding symplectic blow-up 
$\widetilde{X}$ has a nontrivial $n$-tuple Massey product
$\langle\widetilde{S}_1,\dots,\widetilde{S}_n\rangle$, where
$\widetilde{S}_i\in N(H^3(\widetilde{X}))$, $1\leq i\leq n$.}

\smallskip

{\bf Theorem C} {\it Let a symplectic manifold $(Y,\omega)$ have a nontrivial
triple matrix Massey product. Then for any symplectic embedding $Y\subset X$
of codimension greater or equal than $8$ the corresponding symplectic blow
up $\widetilde{X}$ also has a nontrivial triple matrix Massey product.}

\smallskip

A particular case of Theorem C was established by us in 
\cite{BT0,BT} and moreover quite general reasonings applied to the blow-ups of
$\C P^n$ along the embedded Kodaira--Thurston manifolds 
give a bound for the codimension equal to $8$ \cite{BT0}.
Later we showed that particularly
for this manifold the bound for the codimension is decreased to $6$.

In \cite{BT} we posed a problem of inheritance of nontrivial Massey products 
in a general case. It was discussed in \cite{RT} where Theorem C for ordinary
Massey products was proved.

\smallskip

{\bf Theorem D} {\it Let  a symplectic manifold $(Y,\omega)$ have a strictly 
irreducible quadruple matrix Massey product 
$\langle S_1, S_2,S_3,S_4\rangle$. Then for any symplectic submanifold 
$Y\subset X$ such that
$$
\codim Y>2\,\sdeg\langle S_1,S_2,S_3,S_4\rangle
$$
the corresponding symplectic blow-up $\widetilde{X}$ has a nontrivial 
quadruple matrix Massey product.}

\smallskip

Theorems A, B, C, and D enables us to construct new manifolds with
nontrivial Massey products by the symplectic blow-up of $X$ along a
submanifold $Y$. In this event it needs that $X$ or $Y$ have nontrivial 
Massey products. Therefore to make this method work we have to have an initial
family of symplectic manifolds with nontrivial Massey products. 
An existence of such a family is guaranteed by Theorem \ref{th4.2.1}, 
which we shall prove in \S 3:

{\it For any $k$ there exist symplectic manifolds with Massey products whose 
weights are strictly equal to $2k$.}

\smallskip

Examples of such manifolds are constructed in 
Proposition \ref{pr4.2.2}.

These results show that an existence of a symplectic structure does not
strongly restrict the homotopy type of a manifold. We would like to
propose the following

\smallskip

{\bf Conjecture} \ 
{\it For any finite polyhedron $P$ and any natural $N$ there are
a symplectic manifold $X$ and an embedding $f: P \to X$ such that
$$
f_{\ast}: \pi_k(P) \to \pi_k(X) \ \ \mbox{is a monomorphism for $k \leq N$.}
$$ 
}

For $N=1$ and $\dim P =2$
the conjecture follows from Gompf's results \cite{G}.

The authors thank D. van Straten and A. Tralle for helpful discussions.

\begin{center}
{\bf \S 2. Massey products and some constructions of the homotopy theory}
\end{center}

{\bf 2.1. The algebra of forms and its minimal model.}

Let $\A$ be a differential graded algebra over a field $k$. This means that 
$\A$ is a direct sum of the subspaces
$\A^l$ formed by homogeneous elements of degree $l\geq 0$:
$$
\A=\oplus_{l\geq 0}\A^l,\ \ \deg a=l \ \ \mbox{for}\ \ a\in\A^l,
$$
and there are the following linear operations on $A$:
an associative multiplication
$$
\wedge:\A^l\times\A^m\to\A^{l+m},\ \ l,m\geq 0,
$$
and a differential
$$
d:\A^l\to\A^{l+1},\ \ l\geq 0.
$$
Moreover it is assumed that the following conditions hold:

1) $a\wedge b=(-1)^{lm}\,b\wedge a$ for $a\in\A^l$, $b\in\A^m$;

2) $d\,(a\wedge b)=d\,a\wedge b+(-1)^l a\wedge d\,b$ for $a\in\A^l$ 
(the Leibniz rule);

3) $d^2=0$.

A homomorphism of differential graded algebras
$(\A,d_{\A})$ and $(\B,d_{B})$ over a field $k$ is
a $k$-linear mapping
$$
f:\A\to\B,
$$
which respects the grading:
$$
f(\A^l)\subset \B^l\ \ \mbox{for}\ \ l\geq 0,
$$
and satisfies the following conditions:
$$
f(a\wedge b)=f(a)\wedge f(b),\ \ a,b\in\A,
$$
$$
d_{\B}\,f(a)=f(d_{\A}\,a),\ \ a\in\A.
$$
The cohomologies of the algebra
$(\A,d_{\A})$ are defined in the natural manner:
$$
H^l(\A,d_{\A})=\Ker(d_{\A}:\A^l\to\A^{l+1})/\Im(d_{\A}:\A^{l-1}\to\A^l).
$$
To every closed element (or cocycle)  
$a\in\A$, i.e. such that $d_{\A}\,a=0$, there corresponds its 
cohomology class $[a]\in H^+(\A,d_{\A})$. If $[a]=0$, then it is said that
$a$ is exact.

A homomorphism 
$f:(\A,d_{\A})\to(\B,d_{\B})$ maps closed elements to closed ones and 
exact elements to exact ones. Therefore it induces a homomorphism
$$
f^{\ast}:H^{\ast}(\A,d_{\A})\to H^{\ast}(\B,d_{\B})
$$
by the formula $f^{\ast}[a]=[f(a)]$.

An algebra  $(\A,d)$ is connected if $H^0(\A,d)=k$. If in addition
$H^1(\A,d)=0$, then it is said that $A$ is simply connected.

It is also assumed that $\A$ is augmented, which means that there is 
an epimorphism of differential graded algebras
$$
\varepsilon:\A\to k
$$
where $k$ is formed by zero degree elements and is endowed with the zero
differential.
The ideal 
$I=\Ker\varepsilon$ is called the augmentation ideal.

In this paper the main example of such an algebra would be the algebra
$\E^{\ast}(X)$ of smooth differential forms on a smooth manifold
$X$. If $X$ is simply connected then $\E^{\ast}(X)$ is also simply connected.

Let $i:x\to X$ be the embedding of the point 
$x\in X$ into  $X$. The induced mapping defines an augmentation
$$
\varepsilon=i^{\ast}:\E^{\ast}(X)\to\E^{\ast}(x)=k=\R.
$$
The augmentation ideal consists of all forms of positive degree and
smooth functions $\varphi:X\to \R$ such that $\varphi(x)=0$.

In the sequel we assume for simplicity that $k = \R$.

A differential graded algebra $(\M,d_{\M})$ is called minimal if

1) $\M^0=\R$, $d\,(\M^0)=0$, and the multiplication by elements of
$\M^0$ coincides with the multiplication by elements of the main field
$k = \R$;

2) $\M^+=\oplus_{l>0}\M^l$ is freely generated by homogeneous elements
$x_1,\dots$, $x_n,\dots$:
$$
\M^+=\wedge(x_1,\dots),
$$
for any $l>0$ there are finitely many such generators of degree $l$, and
$$
\deg x_i\leq\deg x_j \ \ \mbox{for}\ \ i\leq j;
$$

3) the differential $d$ is reducible:
$$
d\,x_i\in\wedge(x_1,\dots,x_{i-1})\ \ \mbox{for} \ i\geq 1,
$$
i. e., $d\,x_i$ is a polynomial in $x_1,\dots,x_{i-1}$.

It is clear that a minimal algebra $\M$ is simply connected if and only of
$\M^1 = 0$. In this case $\deg x_i\geq 2$ for $i\geq 1$ and the
reducibility condition is written as
$$
d\,(\M^+)\subset \M^+\wedge\M^+
$$
where $\M^+\wedge\M^+$ is a linear subspace generated by reducible elements.

An algebra $(\M,d_{\M})$ is called a minimal model for
an algebra $(\A,d_{\A})$ if 

1) the algebra $(\M,d_{\M})$ is minimal;

2) there is a homomorphism $h:(\M,d_{\M})\to (\A,d_{\A})$, which induces an 
isomorphism of the cohomology rings:
$$
h^{\ast}:H^{\ast}(\M,d_{\M})\to H^{\ast}(\A,d_{\A}).
$$

In the sequel homomorphisms satisfying these conditions are called
{\it quasiisomorphisms}.

The fundamental theorem of Sullivan \cite{S,DGMS} reads

\begin{theorem} $($Sullivan$)$ If $(\A,d_{\A})$ is a simply connected
differential graded algebra such that $\dim H^l(\A)<\infty$ for any
$l\geq 0$, then there is a minimal model for $\A$ which is unique up 
to isomorphism.
\end{theorem}

An example of an algebra satisfying the hypothesis of the theorem is
the algebra of smooth forms on a simply connected compact manifold.

If $\A=\E^{\ast}(X)$, then the minimal model $\M_X$ for $\A$ is 
also called the (real) minimal model for the space $X$. Since, by the de Rham,
theorem the cohomologies of the algebra of smooth forms on $X$ are
isomorphic to the real cohomologies of $X$, the homomorphism
$$
h:(\M_X,d_{\M})\to(\E^{\ast}(X),d_X)
$$
induces the isomorphism
$$
h^{\ast}:H^{\ast}(\M_X)\to H^{\ast}(X)
$$
where $H^{\ast}(X)=H^{\ast}(X,\R)$.

{\sl Examples.} 1) Let $X=S^{2n+1}$, $n\geq 1$. Then $\M^+_X$ is generated by 
an element $x$ with $d\,x=0$ and $\deg x=2n+1$.

2) For $X=S^{2n}$, $n\geq 1$, its minimal model $\M^+_X$ is generated by 
elements $x$ and $y$ with $\deg x=2n$, $\deg y=4n-1$, $d\,x=0$, and
$d\,y=x^2$.

An important property of minimal models is given by the following statement
\cite{S,DGMS}.

\begin{theorem} $($Sullivan$)$ 
Any smooth mapping 
$$
f:X\to Y
$$
induces a homomorphism of the minimal models
$$
\widehat{f}:\M_Y\to\M_X
$$
such that the diagram
$$
\begin{array}{ccc}
H^{\ast}(\M_X) & \stackrel{\widehat{f}^{\ast}}{\longleftarrow} & 
H^{\ast}(\M_Y) \\
h^{\ast}_X\downarrow &  & \downarrow h^{\ast}_Y  \\
H^{\ast}(X) & \stackrel{f^{\ast}}{\longleftarrow} &
H^{\ast}(Y)
\end{array}
$$
is commutative.
\end{theorem}

Other properties of minimal models are exposed in 
\cite{S,DGMS} and we note only one of them: there is the natural homomorphism
$$
\Hom(\pi_l(X),\R)=\M^l_X/\M^+_X\wedge\M^+_X,
$$
i. e., irreducible elements from $\M^l_X$ are in a one-to-one
correspondence with homomorphisms from $\pi_l(X)$ to $\R$. 
Here we assume that there is the minimal model for $X$ and
this model is constructed only for simply connected and nilpotent non-simply
connected spaces.

In the non-simply connected case we confine to a description of the minimal 
models for nilmanifolds which are quotients of nilpotent Lie groups
$G$ with respect to uniform lattices $\Gamma \subset G$.

Let $G$ be a simply connected nilpotent Lie group, let
$\G$ be its Lie algebra, and let  $\G^{\ast}$ be the algebra dual to $\G$,
i. e., the algebra of linear functions $f:\G\to\R$.

Let $e^1,\dots,e^n$ be a basis for $\G$ and let $\omega_1,\dots,\omega_n$ 
be the dual basis for $\G^{\ast}$: $\omega_i(e^j)=\delta^j_i$. The structure
constants $c_k^{ij}$ are defined as follows:
$$
[e^i,e^j]=\sum_{k=1}^n c_k^{ij}e^k
$$
and, by the Mal'tsev theorem \cite{Mal}, 
the group $G$ contains a discrete subgroup 
$\Gamma$ with the compact quotient space $G/\Gamma$ if and only if the
structure constants are rational in a 
certain basis $e^1,\dots,e^n$.

The group $G$ is diffeomorphic to $\R^n$ and its elements may be identified 
with vectors $\xi=\sum_i\,\xi_i e^i\in\R^n$. The multiplication is defined by
the Campbell--Hausdorff formula 
$$
\left(\sum_i\,\xi_i e^i\right)\times\left(\sum_j\,\eta_j e^j\right)=
\sum_k\,P_k(\xi_1,\dots,\xi_n,\eta_1,\dots,\eta_n)e^k
$$
and it follows from the nilpotence of $\G$ that 
$P_1,\dots,P_n$ are polynomials. If $c_k^{ij}$ are rational, then
$P_1,\dots,P_n$ are polynomials with rational coefficients and the subgroup
generated by the vectors $e^1,\dots,e^n$ is a (uniform) lattice, i. e., it is
a discrete subgroup with the compact quotient $G/\Gamma$.

Let $\M_X$ be the minimal model for  $X=G/\Gamma$ where $\Gamma$ is a lattice 
in $G$ with the compact quotient $G/\Gamma$. Then 
$\M^+_X$ is generated by elements  $a_1,\dots,a_n$ where
$$
\deg a_1=\dots=\deg a_n=1,
$$
$$
d\,a_k=\sum_{i,j}\,c_k^{ij}\,a_i\wedge a_j.
$$
The mapping $\M_X\to\E^{\ast}(X)$ has the form
$$
a_k\to\omega_k
$$
where $\omega_k$ is the left-invariant form on $X$ generated by the element 
$\omega_k \in \G^{\ast}$. By the Nomizu theorem \cite{N}, this mapping induces
the isomorphism
$$
H^{\ast}(\M_X)\simeq H^{\ast}(X)
$$
and the minimal model for $X$ coincides with the complex generated by 
left-invariant $1$-forms. 

Notice that the minimal models for nilmanifolds are not simply connected.

{\sl Example.} Let $G$ be the group of all matrices of the form
$$
\left(\begin{array}{ccc}
1 & x & z \\
0 & 1 & y \\
0 & 0 & 1 \end{array}\right),\ \ x,y,z\in\R
$$
and let $\Gamma=G_{\Z}$ be a subgroup of all elements with $x,y,z\in\Z$. 
Take the following basis for leftinvariant $1$-forms:
$$
\omega_1=d\,x,\ \ \omega_2=d\,y,\ \ \omega_3=x\,dy-dz.
$$
The minimal model $\M_X$ for the Heisenberg nilmanifold 
$X=G/G_{\Z}$ is generated by elements 
$a_1,a_2,a_3$ of degree $1$ where
$$
d\,a_1=d\,a_2=0,\ \ d\,a_3=a_1\wedge a_2.
$$

{\bf 2.2. Formality.}

To every minimal algebra $(\M,d)$  there corresponds its cohomology ring 
$\H^{\ast}(\M)$, which we consider as a differential graded algebra 
with the zero differential: $(\H^{\ast}(\M),0)$.

If there is a homomorphism
$$
f:(\M,d)\to (\H^{\ast}(\M),0)
$$
inducing an isomorphism of the cohomology rings
$$
f^{\ast}:\H^{\ast}(\M)\simeq \H^{\ast}(\M),
$$
then the minimal algebra  $\M$  is called {\it formal}.

An existence of such an isomorphism means that the algebra
$(\M,d)$  is the minimal model for its cohomology ring. The construction of 
the minimal model for any algebra $(\A,d_{\A})$ is done effectively by 
induction and formality means that $\M$ is reconstructed from $\H^{\ast}(\M)$.

A differential graded algebra $(\A,d_{\A})$ is called 
{\it formal} if its minimal model $\M(\A)$ is formal. If that holds the 
minimal model for $(\A,d_{\A})$ is reconstructed from 
$\H^{\ast}(\A,d_{\A})=\H^{\ast}(\M(\A),d_{\M})$.

If the algebra $\E^{\ast}(X)$ of smooth forms on a smooth manifold $X$ is 
formal, then the space $X$ is called formal itself.

A notion of formality is generalized for compact simply connected and
nilpotent polyhedra. In this case the algebra  $\E^{\ast}(X)$ is replaced
by the algebra of piecewise polynomial forms \cite{S,GM}.

Examples of formal spaces are compact symmetric spaces and compact K\"ahler
manifolds \cite{DGMS}, and also simply connected compact manifolds (without
boundary) of dimension $\leq 6$ \cite{Mil,NM}.

An example of a nonformal algebra is the minimal model for the 
three-dimensional Heisenberg nilmanifold exposed in \S 2.1.
Indeed, if there is a homomorphism 
$$
f:(\M_X,d)\to(H^{\ast}(X),0)
$$
inducing the isomorphism of the cohomology rings, then 
$$
f(a_1)\neq 0,\ \ f(a_3)=0
$$
and, therefore, $f(a_1\wedge a_3)=0$. But the element  $a_1\wedge a_3$   
realizes a nontrivial cohomology class.

The formality criterion for minimal algebras was given in
\cite{DGMS}. To explain it notice that a minimal algebra is
a tensor product of exterior algebras  
$\Lambda(V_l)$ where $V_l$ are the spaces generated by 
$l$-dimensional generators:
$$
\M=\otimes_{l\geq 0}\Lambda(V_l).
$$
In each space $V_l$ let us take a subspace $C_l$ of all closed elements.

\begin{theorem} $($\mbox{{\rm\cite{DGMS}}}$)$ A minimal algebra
$(\M,d)$ is formal if and only if in each $V_l$ there is a complement  
$N_l$  to $C_l$:
$$
V_l=C_l\oplus N_l
$$
such that any closed element from the ideal
$I_N=I(\oplus N_l)$ generated by $N_1,N_2,\dots$ is exact.
\end{theorem}

It is not always convenient to use this criterion and hence the following 
test is used quite often: to find in $H^{\ast}(\M,d)$ nontrivial classes
represented by Massey products. If there are such classes, then the algebra
$(\M,d)$ is not formal. The converse statement is not true as is shown
by an example from \cite{HS}. We shall discuss relations of Massey products to formality later after defining these operations.

\smallskip

{\bf 2.3. The Maurer--Cartan equations in differential algebras and
generalized Massey products.}

Let $(\A,d)$ be a differential graded algebra over $k(=\R)$ with 
an augmentation. Let $M(\A)$ be a set of all upper triangular half-infinite 
matrices with entries from $\A$, zeroes at the diagonal and finitely many
nonzero entries.

Hence, for any matrix
$$
A=(a_{ij})_{i,j\geq 1}\in M(\A)
$$
we have 
$$
a_{ij}\in\A,\ \ a_{ij}=0\ \mbox{for}  \ \ j\leq i\ \mbox{and}  \ \ i,j\geq n+1
$$
for some $n$. Notice that the condition
$$
a_{ij}=0\ \mbox{for}  \ \ i,j\geq n+1
$$
distinguishes in $M(\A)$  a subset $M_n(\A)$ consisting of all 
$(n\times n)$-matrices with entries from $\A$.
We have a natural filtration  
$$
M_1(\A)\subset M_2(\A)\subset\dots\subset M_n(\A)\subset\dots M(\A).
$$

There are the operations of addition and multiplication endowing 
$M(\A)$ with an algebra structure over $k=\R$. 
The subsets  $M_n(\A)$, $n\geq 1$, are subalgebras of $M(\A)$.

Let us define a differential $d$ on $M(\A)$ as 
\begin{equation}
d\,A=(d\,a_{ij})_{i,j\geq 1}.
\label{1.2.1}
\end{equation}

The mapping acting on homogeneous elements $a \in A^k$ 
as follows 
$$
a\to\bar{a}=(-1)^ka
$$
generates an automorphism of order $2$ (an involution) of $\A$.  
This automorphism is extended to an automorphism of $M(\A)$ as
$\bar{A}=(\bar{a}_{ij})_{i,j\geq 1}$ and the differential   
(\ref{1.2.1}) satisfies the generalized Leibniz rule
$$
d\,(AB)=(d\,A)B+\bar{A}(d\,B).
$$
The algebra $M$ is naturally bigraded:
$$
M=\sum_{p\geq 1,\,k\geq 0} M^{p,k}.
$$
Indeed, let $(a)_{ij}$ be a matrix with one and only one nonzero entry 
such that it is the $(i,j)$ entry which equals
$a \in \A$. Then define $M^{p,k}$ as a subspace of $M$ generated by 
$$
(a)_{i,(i+p)},\ \ \ \ a\in\A^k,\ \ i=1,2,\dots.
$$
The multiplication acts as follows       
$$
M^{p,k}\otimes M^{q,l}\stackrel{\cdot}{\to} M^{p+q,k+l}
$$
and is not graded commutative.

Another structure on $M$ is the bigraded Lie brackets
$$
M^{p,k}\otimes M^{q,l}\stackrel{[\,,]\:}{\longrightarrow} M^{p+q,k+l}
$$
defined on homogeneous elements as    
$$
[A,B]=A\cdot B-(-1)^{kl}B\cdot A,\ \ A\in M^{p,k},\ \ B\in M^{q,l}.
$$
It is straightforwardly checked that these brackets endow $M$ with a structure
of a differential Lie superalgebra with standard properties:

1. $[A,B]=-(-1)^{kl}[B,A]$;

2. $d\,[A,B]=[d\,A,B]+[\bar{A},d\,B]$;

3. $(-1)^{km}[[A,B],C]+(-1)^{lk}[[B,C],A]+(-1)^{lm}[[C,A],B]=0$
for $A\in M^{p,k}$, $B\in M^{q,l}$, $C\in M^{r,m}$.

Now for any (not only for homogeneous) matrix $A\in M$ define its 
{\it kernel}, $\Ker A$, as a certain two-sided $\A$-submodule of 
$M$. By the definition, $\Ker A$ is linearly 
generated as an $\A$-module by the 
matrices $(1)_{ij}$ such that $A\cdot(1)_{ij}=(1)_{ij}\cdot A=0$.
This implies that
$$
AB=BA=0
$$
for any matrix $B\in\Ker A$.

Notice that generically $\Ker A$ is not an ideal with respect to either
the multiplication, either the Lie brackets. But there is a submodule 
$\Ker^\prime A$ of $\Ker A$ which is a two-sided ideal with respect to 
both these operations. Let us describe it. 
Let $m$ be the number of the first nonzero column in $A$ 
and let $k$ be the number of the last nonzero row in $A$. Then define 
$\Ker^\prime A$ as a $\A$-submodule of $M$ generated by the matrices
$(1)_{ij}$ with $i<m$ and $j>k$.

Let us put the Maurer--Cartan operator $\mu: M\to M$ equal to 
$$
\mu(A) = d\,A - \bar{A} \cdot A.
$$

\begin{definition}
\label{def1.2.1a}
A matrix $A\in M$ is called the matrix of a formal connection if it
satisfies the Maurer--Cartan equation in $\A$:
\begin{equation}
d\,A-\bar{A}\cdot A\equiv 0\ \ \mod\Ker A.
\label{1.2.2}
\end{equation}
\end{definition}

In other words, $A$ is a formal connection if $\mu(A) \in \Ker A$.

\begin{definition}
\label{def1.2.1b}
Let $A$ be a formal connection. Then $\mu(A)$ is called the curvature 
matrix of $\A$.
\end{definition}

{\sl Example.}  Let $A$ be an upper triangular $(n+1)\times(n+1)$-matrix 
with zeroes on the diagonal of the form     
\begin{equation}
A=\left(\begin{array}{cccccc}
0 & a_1 & \ast   & \dots & \ast & \ast \\
0 & 0   & a_2 & \dots & \ast & \ast \\
&       &     & \dots &   &   \\
0 & 0   & 0   & \dots & 0 & a_n  \\
0 & 0   & 0   & \dots & 0 & 0
\end{array}\right),\label{matrix}
\end{equation}
where $a_1,\dots,a_n$ do not vanish. Then $\Ker A\cap M_{n+1}(\A)$   
is generated over $k$ by the matrix with one and only one nonzero entry,
which is the $(n+1,n+1)$ entry and equals $1$.
The condition (\ref{1.2.2}) takes the form    
\begin{equation}
d\,A-\bar{A}\cdot A=
\left(\begin{array}{cccc}
0 & \dots & 0 & \tau \\
0 & \dots & 0 & 0 \\
& \dots& & \\
0 & \dots & 0 & 0
\end{array}\right).
\label{star}
\end{equation}

\begin{lemma} 
\label{lem1.2.2} 
$($The Bianci identity$)$
For any matrix $A\in M$ we have  
$$
d\,\mu(A)=\overline{\mu(A)}\cdot A-A\cdot\mu(A).
$$
\end{lemma}

{\sl Proof.}  It follows from the definition of $\mu$ that 
$$
d\,\mu(A)=-d\,(\bar{A}\cdot A)=
-d\,\bar{A}\cdot A-A\cdot d\,A=\overline{d\,A}\cdot A-A\cdot d\,A=
$$
$$
=\overline{(\mu(A)+\bar{A}\cdot A)}\cdot A-A(\mu(A)+\bar{A}\cdot A)=
$$
$$
= \overline{\mu(A)}\cdot A+
A\cdot\bar{A}\cdot A-A\cdot\mu{A}-A\cdot\bar{A}\cdot A=
$$
$$
=\overline{\mu(A)}\cdot A-A\cdot\mu(A).
$$
Lemma is proved.

\begin{corollary} 
\label{cor1.2.3}
The curvature matrix of a formal connection is closed.
\end{corollary}

{\sl Proof.} If $A$ is a formal connection, then $\mu(A)\in\Ker A$, and,
therefore, $\overline{\mu(A)}\in\Ker A$. Now
the Bianci identity implies that $d\,\mu(A)=0$.

Let us denote the cohomology class of a closed element
$a\in\A$, $d\,a=0$, by $[a]\in H^{\ast}(\A)$. 
If $A\in M(\A)$ is a closed matrix, $d\,A=0$, then we denote by  
$[A]=([a_{ij}])_{i,j\geq 1}\in M(H^{\ast}(\A))$ the corresponding matrix 
whose entries are the cohomology classes of $a_{ij}$.

Corollary \ref{cor1.2.3} guarantees the correctedness of 
the following definition.   

\begin{definition} 
\label{def1.2.4}
\footnote{Notice that the Massey products are defined in more general 
situation when $\A$ is an algebra over a commutative ring (for instance, $\Z$)
and the multiplication is associative but not necessary graded 
commutative ($a\wedge b=(-1)^{pq}b\wedge a$). These are, for instance,
the classical Massey products in $H^{\ast}(X,\Z)$.}
Let $A$ be a solution to the Maurer--Cartan equation. Then the entries of the
matrix $[\mu(A)]$ are called generalized Massey products.
\end{definition}

In other words, the generalized Massey products are the cohomology classes
of the curvature matrices of formal connections. Speaking in the 
differential-geometric language, we say that the generalized Massey products
measure the deviations of connections from flat ones.
If all generalized Massey products vanish, then the connection is flat
(see \S 2.4).

{\sl Remarks.}

1) The classical $n$-tuple Massey products (see \S 2.3) differ by sign from 
the corresponding generalized Massey products. We changed the sign in 
Definition \ref{def1.2.1a} to stress that the constructions of abstract 
homological algebra and differential geometry are parallel.

2) It is shown that if $\alpha\in H^\ast(\A)$ is a generalized Massey
product, then $t\alpha$, where $t\in k$, is also a generalized Massey
product. By the equality $[\mu(A)]=[d\,A-\bar{A}\cdot A]=-[\bar{A}\cdot A]$,
we conclude  that the generalized Massey products are the entries of 
the matrix $[\bar{A}\cdot A]$.

From two equivalent notions which are 
``a solution to the Maurer--Cartan equation'' 
and ``a formal connection'' we shall use in the sequel
only the first one.

Let $\varepsilon:\A\to k$ be an augmentation and let $I=\Ker\varepsilon$ be
the augmentation ideal in $\A$. 
For a connected algebra $\A$, i. e. $\A^0=k$, the augmentation ideal 
$I=\A^+$ is an ideal of elements of positive degree and this takes place,
for instance, for a minimal algebra $\A$. By natural reasons of homological
algebra, which are also used at the definition of the reduced bar-construction,
we call a solution to the Maurer--Cartan equation {\it nontrivial} if 
it belongs to the ideal $M(I)$. Otherwise it is called trivial.
In the sequel we shall consider only nontrivial solutions.   

{\sl Examples.}

{\it 1$)$ Solutions of the Heisenberg type and the ideal
$H^+(\A)\cdot H^+(\A)$.}

Let us consider a matrix
\begin{equation}
A=\left(\begin{array}{cccccc}
  0 & \bar{a}_1 & \bar{a}_2 & \dots  & \bar{a}_l & c      \\
  0 & 0         &           &        & 0         & b_1    \\
    &           & \ddots    &        &           & b_2    \\
    &           &           & \ddots &           & \vdots \\
  0 &           & \dots     &        & 0         & b_l     \\
  0 &           & \dots     &        & 0         & 0
  \end{array}\right)
\label{1.2.3}
\end{equation}
where $a_i,b_j,c\in I$, $a_{ij}=0$, $i,j>l+1$. We say that this matrix
is of the Heisenberg type keeping in mind an analogy with  
elements of the Lie algebra of the generalized Heisenberg group. Let us 
try to find a solution of the form (\ref{1.2.3}) to the Maurer--Cartan
equation. It follows from the definition that 
$d\,a_i=d\,b_j=0$ and if  $\alpha_i$, $\beta_j$ are the corresponding 
cohomology classes, then the generalized Massey product defined by
$A$ is
$$
[\mu(A)]=-\sum^n_{i=1}\alpha_i\beta_i\in H^+(\A)\cdot H^+(\A).
$$
Hence we see that any element of the ideal $H^+(\A)\cdot H^+(\A)$ 
is represented by a Massey product of the Heisenberg type.
In the classical interpretation (see below) this is a double matrix
product
$$
\langle\:(\bar{\alpha}_1,\dots,\bar{\alpha}_n),\left(\begin{array}{c}
                                             \beta_1 \\
                                             \vdots  \\
                                             \beta_n
                                             \end{array}\right)\rangle.
$$

\begin{definition} 
\label{def1.2.5}
Generalized Massey products belonging to  
$H^+(\A)\cdot H^+(\A)$ are called completely reducible. A generalized
Massey product is called completely irreducible if it 
has a nontrivial image under the projection
$$
H^\ast(\A)\to H^\ast(\A)/H^+(\A)\cdot H^+(\A).
$$
\end{definition}

{\sl Remark.} We shall see below that there are relations between 
Massey products which sometimes enables us to reduce complex products to more
simple ones and, in particular, to ordinary products. As well as the notion of 
nontriviality of Massey products, which will be introduced below, 
the complete irreducibility characterizes a ``nonsimplifiable'' property of
Massey products.

{\it 2$)$ Triple $($ordinary$)$ Massey products.}

Let $\alpha$, $\beta$, and  $\gamma$ be the cohomology classes of closed 
elements $a \in\A^p$, $b\in\A^q$, and  $c\in\A^r$. The triple Massey product   
$\langle\alpha,\beta,\gamma\rangle$ is defined if 
the Maurer--Cartan equation is
solvable for
$$
A=\left(\begin{array}{cccc}
0 & a & f & h \\
0 & 0 & b & g \\
0 & 0 & 0 & c \\
0 & 0 & 0 & 0
\end{array}\right).
$$
This equation is equivalent to  
\begin{equation}
d\,f=(-1)^p\,a\wedge b,\ \ d\,g=(-1)^q\,b\wedge c
\label{ast}
\end{equation}
and that implies that the product $\langle\alpha,\beta,\gamma\rangle$    
is defined if and only if  
$$
\alpha\cup\beta=\beta\cup\gamma=0\ \ \mbox{in}\ \ H^{\ast}(\A).
$$
The matrix $\mu(A)$ has the form (\ref{star}) for  $n=3$ and defines
the Massey product $[\mu(A)]$ which equals  
$$
\langle\alpha,\beta,\gamma\rangle=[\tau]=
[(-1)^{p+1} a\wedge g+(-1)^{p+q} f\wedge c].
$$
Since $f$ and  $g$ are defined by (\ref{ast}) up to closed elements
from $\A^{p+q}$  and $\A^{q+r}$, the triple Massey product 
$\langle\alpha,\beta,\gamma\rangle$ is defined modulo
$\alpha\cdot H^{q+r}(\A)+\gamma\cdot H^{p+q}(\A)$.

{\it 3$)$ $n$-tuple ordinary Massey products} \cite{K}.

Given $\alpha_1,\dots,\alpha_n\in H^{\ast}(\A)$, the $n$-tuple Massey
product $\langle\alpha_1,\dots,\alpha_n\rangle$ is defined if there are 
their representatives $a_1,\dots,a_n$ such that the Maurer--Cartan equation 
is solvable for $A$  of the form (\ref{matrix}). In this event $\mu(A)$
takes the form (\ref{star}) and the $n$-tuple (ordinary) Massey product is 
$\langle\alpha_1,\dots,\alpha_n\rangle=[\tau]$. It is evident
that for $n=3$ we obtain a triple Massey product, which was already introduced
and discussed in detail with explaining the type of multivaluedness.

{\it Initial data.}

Let $A$ be a solution to the Maurer--Cartan equation. By the definition, 
the first diagonal of $A$ consist only of closed elements, i.e.,
$d\,a_{i,i+1}=0,\ \ i=1,2,\dots$. Other diagonals may also contain closed 
elements. This happens, for instance, for matrix Massey products which we
shall discuss below. Consider the classes of all closed elements of $A$ or 
a part of them and call them the initial data of the corresponding 
Maurer--Cartan problem.

The following statement, which is well-known for classical $n$-tuple ordinary
or Massey products, shows that the set of the cohomology classes 
$[\mu(A)]$ of different solutions with the given initial data depends only
on the initial data.          

\begin{proposition} 
\label{pr1.2.5}
Let  $A$ be a solution to the Maurer--Cartan equation and let 
$\K=\Ker A$. Any replacement of the entry 
$a_{ij}$  of the matrix $A$ by $a_{ij}+d\,b$
is completed to a replacement of $A$ by $A^\prime$ 
such that 

$1$. $A^\prime=A+(d\,b)_{ij}\ \ \mod M^{j-i+1,\ast}$;

$2$. $A^\prime$ is a solution to the Maurer--Cartan equation $\mod\K$;

$3$. $[\mu(A^\prime)]=[\mu(A)]$.
\end{proposition}

{\sl Proof.}  A straightforward computation shows that the conditions 
1 and 2 are satisfied by
$$
A^\prime=A+(d\,b)_{ij}+A\cdot(b)_{ij}-(\bar{b})_{ij}\cdot A.
$$
It is also computed that
$$
\mu(A^\prime)=\mu(A)+d\,((A\cdot(b)_{ij}-(\bar{b})_{ij}\cdot A)\cap\K)
$$
and, therefore, the condition 3 also holds for this choice of $A^{\prime}$.
This proves the proposition.

{\it Naturality of generalized Massey products.}

Let  $f:\A\to\B$ be a homomorphism of differential graded algebras.
It induces a mapping  $\widehat{f}: M(\A)\to M(\B)$ by the formula  
$$
\widehat{f}((a_{ij})_{i,j\geq 1})=(f(a_{ij}))_{i,j\geq 1}.
$$
It is clear that for any matrix $A\in M(\A)$ we have
$$
\widehat{f}(\Ker A)\subseteq\Ker(\widehat{f}(A)).
$$
Therefore, any solution $A$ to the Maurer--Cartan equation in 
$\A$ is mapped into a solution 
$\widehat{f}(A)$  of the Maurer--Cartan equation in $\B$. 
Hence we obtain a mapping of generalized Massey products:
\begin{equation}
f^\ast([\mu(A)])=[\mu(\widehat{f}(A))].
\label{1.2.4}
\end{equation}
Standard reasonings of homological algebra
prove the following statement.

\begin{proposition} 
\label{pr1.2.6}
If $f$ is a quasiisomorphism of differential graded algebras, then
the mapping $(\ref{1.2.4})$ is one-to-one on generalized Massey products. 
\end{proposition}

For matrix Massey products this statement is proved in 
\cite{JPM} (Theorem $1.5$).

This implies that for the minimal model $\M=\M(\A)$ for 
$\A$ the generalized Massey products in $H^{\ast}(\A)$ may be 
computed from the minimal model.

{\bf 2.4. The classical operations: $n$-tuple ordinary and matrix 
Massey products.}

In this subsection we show how the classical Massey products introduced 
in the 60s in \cite{K} and \cite{JPM} fit into the general picture.
It needs to mention that the relation between Massey products and the
Maurer--Cartan equation was first noticed by May \cite{JPM} but this analogy
was not developed.

Let $\alpha_1,\dots,\alpha_n$, $\alpha_i\in H^{p_i}(\A)$, $i=1,\dots,n$,
be cohomology classes and let $a_i\in\alpha_i$ be cocycles representing them in
$\A$.

\begin{definition}
\label{def1.3.1}
A set $A=(a(i,j))$, $1\leq i\leq j\leq n$, $(i,j)\neq(1,n)$, consisting of
some elements of $\A$ is called a defining system for 
$\langle a_1,\dots,a_n\rangle$, if

$1$. $a(i,i)=a_i$, $i=1,\dots,n$;

$2$. $a(i,j)\in\A^{p(i,j)+1}$,  $p(i,j)=\sum^j_{r=i}(p_r-1)$;

$3$. $d\,a(i,j)=\sum_{r=i}^{j-1}\bar{a}(i,r)\cdot a(r+1,j)$.

If these conditions hold, then

$4$. the $(p(1,n)+2)$-dimensional cocycle 
$$
c(A)=\sum_{r=1}^{n-1}\bar{a}(1,r)a(r+1,n)
$$
is called the cocycle of a defining system $A$.
\end{definition}

\begin{definition}
\label{def1.3.2}
An $n$-tuple product $\langle a_1,\dots,a_n\rangle$ is defined if
there is at least one defining system for it. If this product is defined, then
$\langle\alpha_1,\dots,\alpha_n\rangle$ consists of the cohomology
classes $[c(A)]\in H^{p(1,k)+2}(\A)$ of all defining systems $A$ for      
$\langle a_1,\dots,a_n\rangle$.
\end{definition}

In these form these definitions appeared in \cite{K}. We show they relate
to the definitions given in the preceding subsection.     

Let $A$ be a defining system for $\langle a_1,\dots,a_n\rangle$.
We write it as a matrix denoted by the same letter:    
$$
A=\left(\begin{array}{cccccc}
  0      & a(1,1) & a(1,2) & \dots & a(1,n-1)   & a(1,n)   \\
  0      & 0      & a(2,2) & \dots & a(2,n-1)   & a(2,n)   \\
         &        &        & \dots &            &           \\
  0      & 0      & 0      & \dots & a(n-1,n-1) & a(n-1,n) \\
  0      & 0      & 0      & \dots & 0          & a(n,n)   \\
  0      & 0      & 0      & \dots & 0          & 0
  \end{array}\right)
$$
where $a(1,n)$ is an arbitrary element of $\A^{p(1,n)+2}$.
It is easily seen that the conditions 1--3 of Definition \ref{def1.3.1} are
equivalent to a condition that $A$ satisfies the Maurer--Cartan equation.
In this event the matrix $-\mu(A)$ has a single nonzero entry, which equals
$c(A)$. The converse is also true: if $A$ is a solution to the 
Maurer--Cartan equation with the initial data 
$a_{i\,,i+1}=a_i$, $i=1,2,\dots,n$, then the entries of $A$ give a 
defining system for      
$\langle a_1,\dots,a_n\rangle$.

Matrix Massey products were first introduced by May \cite{JPM}. Since the 
matrix entries have to stay homogeneous, we need to define {\it 
multipliable} matrices. 

\begin{definition}
\label{def1.3.3}
Let $N(I)$ be the set of all finite rectangular matrices whose entries are 
homogeneous elements from the augmentation ideal $I$. Two matrices
$X,Y\in N(I)$ of sizes $p\times q$ and $r\times s$ are multipliable
(in this order) if

$1$. $q=r$;

$2$. $\deg x_{ij}+\deg y_{jk}$ is independent of $j$ for all  
$(i,k)$, $1\leq i\leq p$, $1\leq k\leq s$.
\end{definition}

To any matrix $X=(x_{ij})^{\stackrel{p}{q}}_{\stackrel{i=1}{j=1}}\in N(I)$ 
there corresponds an integer matrix  $\D(X)=(\deg x_{ij})$. If $X$  and $Y$
are multipliable, then
$$
\D(X\cdot Y)=\D(X)\ast\D(Y)=(\deg x_{ij}+\deg y_{jk})_{ik}.
$$
The rest part of the definition of matrix Massey products is completely 
analogous to the definition of ordinary products.
Let $V_1$, $V_2$, $\dots$, $V_n$ be matrices from
$ N(H^+(\A))$ such that $V_i$ and $V_{i+1}$ are multipliable for
$i=1,\dots,n-1$. Consider matrices 
$X_i\in N(I)$ consisting of cocycles representing matrices $V_i$, 
$i=1,\dots,n$, i. e.,
$$
d\,X_i=0\ \ \mbox{}\ \ [X_i]=V_i,\ \ i=1,\dots,n.
$$
A sequence $X_i$, $i=1,\dots,n$, is multipliable by its construction and
$\D(X_i)=\D(V_i)$, $i=1,\dots,n$. Let us denote by $K+m$ the matrix obtained 
from the integer matrix $K=(k_{ij})$ by adding to each its entry an 
integer number $m$: $K+m=(k_{ij}+m)$.

\begin{definition}
\label{def1.3.4}
A subset $A=(X(i,j))$, $1\leq i,j\leq n$, $(i,j)\neq (1,n)$, of  $ N(I)$ 
is called a defining system for $\langle X_1,\dots,X_n\rangle$ if  

$1.$ $X(i,i)=X_i$, $i=1,\dots,n$;

$2.$ $\D(X(i,j))=(\D(X_i)-1)\ast(\D(X_{i+1})-1)\ast\dots\ast(\D(X_j)-1)+1$;

$3.$ $d\,X(i,j)=\sum_{r=i}^{j-1}\bar{X}(i,r)\cdot X(r+1,j)$.

If these condition hold, then

$4.$ the matrix cocycle  
$$
c(A)=\sum_{r=1}^{n-1}\bar{X}(1,r)\cdot X(r+1,n)
$$
of degree $\D(1,p)=\D(X_1)\ast\dots\ast\D(X_n)-n+2$ is called the cocycle of 
a defining system $A$.
\end{definition}

\begin{definition}
\label{def1.3.5}
An $n$-tuple matrix product $\langle X_1,\dots,X_n\rangle$ is defined if 
there is at least one defining system for it. If this product is defined, then
$\langle V_1,\dots,V_n\rangle$ consists of the matrix cohomology classes
$[c(A)]$ of all defining systems $A$ for $\langle X_1,\dots,X_n\rangle$.
\end{definition}
As in the case of ordinary products to every defining system $A$ there
corresponds a block matrix
$$
A=\left(\begin{array}{cccccc}
  0 & X(1,1) & X(1,2) & \dots & X(1,n-1)   & X(1,n)   \\
  0 & 0      & X(2,2) & \dots & X(2,n-1)   & X(2,n)   \\
    &        &        & \dots &            &           \\
  0 & 0      & 0      & \dots & X(n-1,n-1) & X(n-1,n) \\
  0 & 0      & 0      & \dots & 0          & X(n,n)   \\
  0 & 0      & 0      & \dots & 0          & 0
  \end{array}\right),
$$
where $X(i,j)$ are blocks described by Definition \ref{def1.3.4} and 
$X(1,n)$ is an arbitrary element of  $\Ker A$. 
The matrix $A$ satisfies the Maurer--Cartan equation and $-\mu(A)$ has 
one and only one nonzero block entry. It is the $(1,n)$ entry, which equals
$c(A)-d\,X(1,n)$:
$$
-[\mu(A)]=[c(A)].
$$
The inverse is also valid: if $A$ is a solution to the Maurer--Cartan 
equation whose initial data are block matrices realizing the matrix cohomology
classes $V_1,\dots,V_n$, then $-[\mu(A)] \in \langle V_1,
\dots,V_n\rangle$.

An analogue of Proposition \ref{pr1.2.5} for $n$-tuple products was proved
in \cite{K} for ordinary products and in \cite{JPM} for matrix products.

\begin{proposition}
\label{pr1.3.6}
Let an $n$-tuple matrix Massey product
$\langle V_1,\dots,V_n\rangle$ be defined. Then its values depends only
on matrix cohomology classes $V_1,\dots,V_n$.
\end{proposition}

{\it Indeterminacy and nontriviality of Massey products.}

Matrix Massey products are not a formal generalization of ordinary
$n$-tuple operations and appear quite naturally, for instance, 
at describing the indeterminacy of an $n$-tuple product. This effect is
almost nonvisible for triple products  $\langle v_1,v_2,v_3\rangle$ but
it is already spectacular for quadruple products.

\begin{definition}
\label{def1.3.7}
Let an  $n$-tuple matrix Massey product $\langle V_1,\dots,V_n\rangle$ 
is defined. Then it is called trivial if $0\in\langle V_1,\dots,V_n\rangle$.
\end{definition}

In accord with Definition \ref{def1.2.5} we give the following  

\begin{definition}
\label{def1.3.7b}
An $n$-tuple $($matrix$)$  Massey product $\langle V_1,V_2,\dots,V_n\rangle$
is called completely reducible if
$$
\langle V_1,\dots,V_n\rangle\subset H^+(\A)\cdot H^+(\A),
$$
otherwise it is called irreducible. Such a product is called
strictly irreducible if
$$
\langle V_1,\dots,V_n\rangle\cap H^+(\A)\cdot H^+(\A)=\emptyset.
$$
\end{definition}

Now assume that $\langle V_1,\dots,V_n\rangle$ and 
$\langle W_1,\dots,W_n\rangle$ are defined and have equal multidegrees,
i. e., $\D(V_1,\dots,V_n)=\D(W_1,\dots,W_n)$.  
Put
$$
\langle V_1,\dots,V_n\rangle+\langle W_1,\dots,W_n\rangle=
$$
$$
=
\{x+y;\: x\in\langle V_1,\dots,V_n\rangle,\: 
y\in\langle W_1,\dots,W_n\rangle\},
$$
$$
\lambda\langle V_1,\dots,V_n\rangle=\{\lambda x;\: 
x\in\langle V_1,\dots,V_n\rangle\}.
$$

\begin{definition}
\label{def1.3.8}
The set
$$
\In\langle V_1,\dots,V_n\rangle=\{x-y;\:x,y\in\langle V_1,\dots,V_n\rangle\}
$$
is called the indeterminacy of the Massey product
$\langle V_1,\dots,V_n\rangle$.
\end{definition}

Hence, if $\langle V_1,\dots,V_n\rangle$ is trivial, then
$$
\langle V_1,\dots,V_n\rangle\subset\In\langle V_1,\dots,V_n\rangle.
$$

The following statement proved in \cite{JPM} describes the indeterminacy of
$n$-tuple products in terms of $(n-1)$-tuple products.

\begin{proposition}
\label{pr1.3.9}
Let a product $\langle V_1,\dots,V_n\rangle$ be defined. Then
$$
\In\langle V_1,\dots,V_n\rangle\subset
\cup_{(X_1,\dots,X_{n-1})}\langle W_1,\dots,W_{n-1}\rangle,
\ \ \mbox{where}
$$
$$
W_1=\left(V_1\ \ X_1\right),\ \
W_k=\left(\begin{array}{cc}
    V_k & X_k \\
    0   & V_{k+1}
    \end{array}\right),\ \ 2\leq k\leq n-2,
$$
$$
W_{n-1}=\left(\begin{array}{c}
        X_{n-1} \\
        V_n
        \end{array}\right).
$$
The union is taken over all $(n-1)$-tuples $X_1,\dots,X_{n-1} \in
N(H^\ast(\A))$ such that
$\D(X_k)=\D(V_k,V_{k+1})-1$. Moreover, for $n=3$ we have   
$$
\In\langle V_1,V_2,V_3\rangle=\cup_{(X_1,X_2)}\langle W_1,W_2\rangle=
\cup_{(X_1,X_2)}(\bar{V}_1\cdot X_2+\bar{X}_1\cdot V_3).
$$
\end{proposition}

{\sl Remark.} 
The properties of Massey products to be nontrivial and to be irreducible are
independent. A Massey product may be completely reducible but not trivial
(see \S 2.3) and if $n>3$, then $n$-tuple Massey products may be trivial but 
not completely irreducible. If $n=3$, then triviality implies reducibility.

\newpage

{\it Strictly defined Massey products.}

To obtain all elements of an $n$-tuple ordinary or matrix product it needs
to consider all defining systems for this product or, which is
equivalent, to consider all solutions to the Maurer--Cartan equation
with given initial data. Defining systems are inductively constructed and 
such a process leads to multivaluedness of the product. 
Moreover it a situation may appear when partially constructed 
defining system can not be completed. In difference 
with a multivaluedness the last effect is not well-known because it takes 
place not for triple but for $n$-tuple products with $n \geq 4$. A potential
impossibility of completing partially constructed defining systems gives rise
to the following definition, which first appeared in \cite{JPM2} (see also
\cite{JPM}).

\begin{definition}
\label{def1.3.10}
Let a matrix product $\langle V_1,\dots,V_n\rangle$ be defined.
Then it is called strictly defined if any product
$$
\langle V_i,\dots,V_j\rangle,\ \ 1\leq j-i\leq n-2
$$
contains only the zero matrix.
\end{definition}

It is clear that triple products are strictly defined.
For many reasons strictly defined products are more useful for
applications. The proof of the following simple lemma may be found in
\cite{JPM}.

\begin{lemma}
\label{lem1.3.11}
A matrix product $\langle V_1,\dots,V_n\rangle$ is strictly defined if and 
only if any its partially determined defining system
$\{X(i,j);\: j-i\leq k\}$, $1\leq k\leq n-2$, 
can be extended to a defining 
system $A$ for $\langle V_1,\dots,V_n\rangle$.
\end{lemma}

{\bf 2.5. Massey products as obstructions to formality.}

The most important for this paper property of Massey operations
is given by

\begin{theorem}
Let a differential graded algebra $\A$ is formal. Then

$a)$ generalized Massey products in $H^\ast(\A)$ are completely 
reducible; 

$b)$ $n$-tuple Massey products in $H^\ast(\A)$ are trivial for $n\geq 3$.
\label{th1.4.1}
\end{theorem}

{\sl Proof.} Let $\M=(\M,d)$ be the minimal model for $\A$ and let 
$H^\ast=H^\ast(\A,d)$. 
By the formality of $A$, there is a quasiisomorphism    
\begin{equation}
h:(\M,d)\to(H^\ast,0).
\label{1.4.1}
\end{equation}
By Proposition \ref{pr1.2.6}, we can compute the generalized Massey products
in $H^{\ast}(\A)$ only from $(\M,d)$. 

a) Let $A$ be a solution to the Maurer--Cartan equation and 
$[\mu(A)]$ be the corresponding   Massey product. By (\ref{1.2.4}), 
we obtain
$$
h^\ast([\mu(A)])=[\mu(\widehat{h}(A))]=
-[\widehat{h}(\bar{A})\cdot\widehat{h}(A)]\in M(H^+\cdot H^+).
$$
Since $h^\ast$ is an isomorphism, we have $[\mu(A)]\in \M(H^+\cdot H^+)$
which proves a).

b) Let $\langle V_1,\dots,V_n\rangle$ be an $n$-tuple matrix  Massey
product. It follows from Proposition \ref{pr1.2.6}  
that
\begin{equation}
\langle V_1,\dots,V_n\rangle=\langle V_1,\dots,V_n\rangle_H
\label{1.4.2}
\end{equation}
where  the symbol $\langle\ \ \rangle_H$ denotes a Massey product
in the differential graded algebra $(H^\ast,0)$.  Since $d\equiv 0$,   ,   
it is clear that the family of matrices
$A=(X(i,i)=V_i;\: X(i,j)=0;\: 1\leq i,j\leq n,\: j-i>0: (i,j)\neq(1,n))$  
is a defining system for  
$\langle V_1,\dots,V_n\rangle_H$. For $n\geq 3$ we have $c(A)=0$ and, 
therefore, $0\in\langle V_1,\dots,V_n\rangle_H$. The equality  
(\ref{1.4.2}) proves b) and hence finishes the proof of the theorem.

\begin{corollary}
\label{cor1.4.2}
If a differential graded algebra $(\A,d)$ has irreducible or
nontrivial matrix Massey products, then the algebra is nonformal.     
\end{corollary}

This reasoning enables us to construct obstructions to formality.
         
{\sl Example.} Let us consider the generalized Heisenberg group,
which is the group of all matrices of the form
$$
\left(\begin{array}{cccccc}
1      & x_1 & x_2    & \dots  & x_n & z      \\
0      & 1   & 0      & \dots  & 0   & y_1    \\
       &     & 1      &        &     & y_2    \\
\vdots &     &        & \ddots &     & \vdots  \\
       &     &        &        & 1   & y_n     \\
0      &     & \dots  &        & 0   & 1
\end{array}\right)
$$
where $x_i,y_j,z\in\R$, $1\leq i,j\leq n$.     
Take for a uniform lattice the subgroup of integer matrices and denote the
corresponding nilmanifold by $X_n$. The minimal model $\M_n=\M(X_n)$ for 
$X_n$ coincides with the complex generated by left-invariant $1$-forms
\cite{N}. 
There is a basis for such forms, which in the coordinates $(\bar{x},\bar{y},z)$
is written as
$$
\alpha^+_i=d\,x_i,\ \ \alpha^-_i=d\,y_i,\ \ 1\leq i\leq n;\  \
\beta=-d\,z+\sum_{i=1}^nx_i\,d\,y_i.
$$
Therefore, we have
$$
\M_n=\Lambda(\alpha^\pm_i,\: 1\leq i\leq n,\: \beta;\; d\,\alpha^\pm_i=0,\:
d\,\beta=\alpha^+_1\wedge\alpha^-_1+\dots+\alpha^+_n\wedge\alpha^-_n).
$$
Using this representation one can compute $H^\ast(X_n,\R)$, but we only notice
that the first nontrivial Massey products appear at the dimension $n+1$ and
are of the form 
\begin{equation}
[\beta\wedge\alpha^{\varepsilon_1}_1\wedge\alpha^{\varepsilon_2}_2
\wedge\dots\wedge
\alpha^{\varepsilon_n}_n]
\label{1.4.3}
\end{equation}
where $\varepsilon_i\in\{\pm 1\}$, $i=1,\dots,n$.
The elements (\ref{1.4.3}) do not lie in $H^+\cdot H^+$ and it is
clear that
$$
[\beta\wedge\alpha^{\varepsilon_1}_1\wedge\dots\wedge
\alpha^{\varepsilon_n}_n]=
$$
$$
=\langle(\varepsilon_1\alpha^{-\varepsilon_1}_1,
\varepsilon_2\alpha^{-\varepsilon_2}_2,\dots,
\varepsilon_n\alpha^{-\varepsilon_n}_n),
\left(\begin{array}{c}
\alpha^{\varepsilon_1}_1 \\
\alpha^{\varepsilon_2}_2 \\
\dots \\
\alpha^{\varepsilon_n}_n
\end{array}\right),
\alpha^{\varepsilon_1}_1\wedge\alpha^{\varepsilon_2}_2\wedge\dots\wedge
\alpha^{\varepsilon_3}_3\rangle.
$$
Hence, the cohomology classes (\ref{1.4.3}) are represented by nontrivial and, 
therefore, irreducible triple matrix Massey products.

{\bf 2.6. Spheric cocycles and the suspension homomorphism.}

Let $X$  be a finite $CW$-complex. An element 
$q\in H^k(X,\R)$ is called  {\it spheric} if there is a mapping 
\begin{equation}
f:S^k\to X
\label{1.5.1}
\end{equation}
such that  $f^\ast(q)\neq 0$.

Since $f^\ast(H^+\cdot H^+)=0$ for any mapping (\ref{1.5.1}), 
the image of a spheric class under the projection      
$$
H^+\stackrel{r}{\to} H^+/H^+\cdot H^+.
$$
does not equal zero.
Let $h_{\ast}:\pi_{\ast}(X)\to H_{\ast}(X,\R)$ be the real Hurewicz 
homomorphism. There is a commutative diagram

\begin{equation}
\label{1.5.2}
\includegraphics{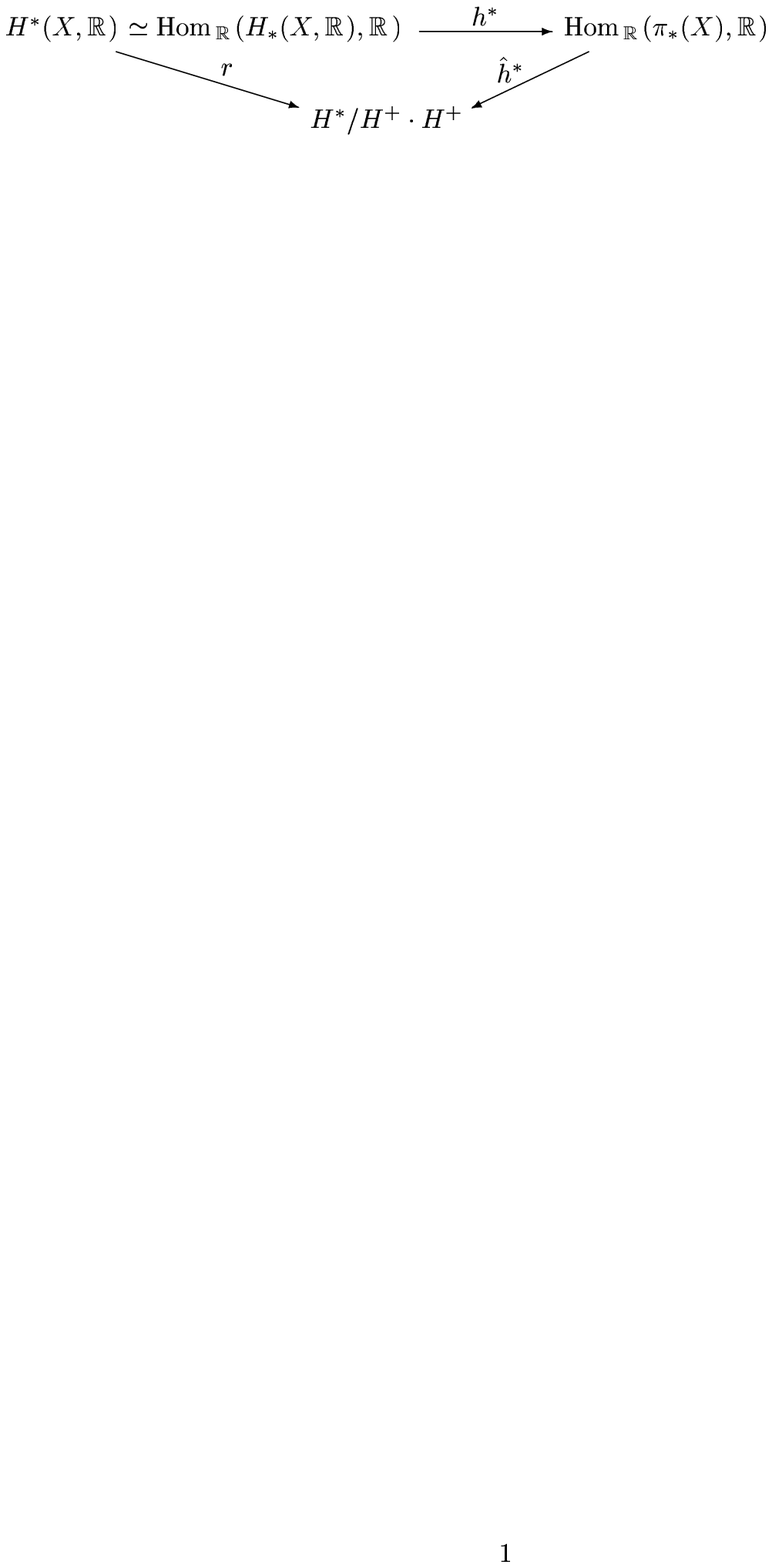}
\end{equation}

Take a subspace $S \subset H^\ast(X,\R)$ complement to $\Ker h$
such that $r\big|_S$ is a monomorphism. Any such a subspace $S$
is called a {\it subspace of spheric generators} in $H^{\ast}(X,\R)$.

It appears a natural question: ``When $H^\ast(X)$
is generated by spheric generators ?'' It follows from the diagram
(\ref{1.5.2}) that that holds if and only if    
$\rk h_{\ast}=\dim(H^\ast/H^+\cdot H^+)$. For fields of the zero 
characteristic this problem is well studied.  Different conditions, which
are equivalent to the condition that the cohomology ring is generated
by spheric generators, may be found in \cite{HS}. 
We shall mention only the following one.

\begin{proposition}
\label{pr1.5.1}
If  $X$ is formal, then $H^\ast(X)$ is generated by spheric generators.
\end{proposition}

The proof of this proposition is also given in \cite{HS}.
Notice that the inverse to this statement is not true. From this proposition
we derive the following obstruction to formality.

\begin{corollary}
\label{cor1.5.2}
If 
$$
\rk h_{\ast}<\dim(H^\ast/H^+\cdot H^+)
$$
where $h_{\ast}$ is the Hurewicz homomorphism,
then $X$ is nonformal.
\end{corollary}

A subspace of spheric generators $S$ and a subspace generated by all
Massey products including ordinary ones are complementary
subspaces of $H^\ast(X,\R)$. To make this statement rigorous let us define
the cohomology suspension homomorphism
\begin{equation}
\Omega^\ast:H^\ast(X,\R)\to H^{\ast-1}(\Omega X,\R)
\label{1.5.3}
\end{equation}
where $\Omega X$ is the loop space of $X$. Denote by $[X,Y]$ 
the set of homotopy classes of mappings from $X$ to $Y$ and
define (\ref{1.5.3}) as follows:
$$
H^n(X,\R)\simeq[X,K(\R,n)]\stackrel{\Omega}{\to}[\Omega X,\Omega K(\R,n)]=
[\Omega X,K(\R,n-1)]\simeq
$$
$$
\simeq H^{n-1}(\Omega X,\R).
$$

An alternative definition comes from considering 
the fibration $p:EX\stackrel{\Omega X}{\to} X$ and the exact cohomology
sequence for $(EX,\Omega X)$. There are homomorphisms 
$$
H^{\ast-1}(\Omega X)\stackrel{\partial}{\to}H^\ast(EX,\Omega X)
\stackrel{p^\ast}{\leftarrow}
H^\ast(X,pt)
$$
and, since $EX$ is contractible, $\partial$ is an isomorphism.
Hence, the homomorphism
$$
\partial^{-1}\circ p^\ast:H^{\ast}(X,pt)\to H^{\ast-1}(\Omega X)
$$
is correctly defined and 
$\Omega^\ast=\partial^{-1}\circ p^\ast$.

There is a commutative diagram   
$$
\begin{array}{ccc}
H^\ast(X,\R) & \stackrel{\Omega^\ast}{\longrightarrow} & 
H^{\ast-1}(\Omega X,\R) \\
h^\ast\downarrow &  & \downarrow h^\ast  \\
\Hom_{\Z}(\pi_{\ast}(X),\R) & \stackrel{\delta^\ast}{\longrightarrow} &
\Hom_{\Z}(\pi_{\ast-1}(\Omega X),\R)
\end{array}
$$
where $\delta^\ast$ is induced by the isomorphism  
$\partial^{-1}:\pi_{k-1}(\Omega X)\to\pi_k(X)$ from the exact homotopy 
sequence for the fibration $EX\stackrel{\Omega X}{\to}X$.  
We derive from this diagram that
$$
\Ker h^{\ast} =\Ker\Omega^{\ast}.
$$
Therefore, there is a decomposition
$$
H^+(X,\R)=S\oplus\Ker\Omega^\ast.
$$

The following statement by May is considered by specialists
to be proved although the proof was never published.

\begin{theorem}
\label{th1.5.3}
The subspace $\Ker\Omega^\ast\subset H^+(X,\R)$   
is generated by all matrix Massey products including ordinary
products.
\end{theorem}

It implies the following useful  

\begin{corollary}
\label{cor1.5.4}
Given a space $X$,  
$$
\Ker h_{\ast}<\dim(H^+/H^+\cdot H^+)
$$
if and only if there are irreducible matrix Massey products in 
$H^\ast(X,\R)$.
\end{corollary}

{\bf 2.7. Massey products and the Eilenberg--Moore spectral sequence.}

Here we expose in brief and in adapted for our study form 
the basic constructions
related to the Eilenberg--Moore spectral sequence. At the end we shall show
a relation between the cohomology suspension homomorphism and this sequence.
Complete algebraic details of that may be found in the paper by Smith
\cite{LS} (see also \cite{HS}).

Let $I$ be the augmentation ideal of a differential graded algebra 
$(\A,\,d)$. The reduced bar-construction for $\A$ is defined as follows.
Put
$$
B\A=\sum_{n\geq 0}\otimes^n I
$$
where the tensor product is considered over the ground field $k$.
If $a_i\in I$ is a homogeneous element of degree $p_i$ for $i=1,\dots,n$,  
then the product $a_1\otimes a_2\otimes\dots\otimes a_n$ is 
denoted by $[a_1|a_2|\dots|a_n]$.
The module $B\A$ has a natural bigrading
$$
\bideg[a_1|a_2|\dots|a_n]=(-n,\sum^n_{i=1}p_i)
$$
and the total degree is defined as
$$
\tdeg[a_1|a_2|\dots|a_n]=\sum^n_{i=1}p_i-n.
$$

By the construction, the module $B\A$ is a Hopf algebra with the diagonal
$$
\Delta([a_1|a_2|\dots|a_n])=
\sum^n_{p=0}[a_1|\dots|a_p]\otimes[a_{p+1}|\dots|a_n]
$$
and the multiplication 
$$
[a_1|\dots|a_k]\cdot[a_{k+1}|\dots|a_{k+n}]=
\sum_{\sigma}\varepsilon_{\sigma}
[a_{\sigma(1)}|a_{\sigma(2)}|\dots|a_{\sigma(k+n)}]
$$
where the summation is taken over all ordered permutations 
$(1,\dots,k+n)$, i. e.,
$\sigma(1)<\sigma(2)<\dots<\sigma(k)$ and 
$\sigma(k+1)<\sigma(k+2)<\dots<\sigma(k+n)$. 
The signature $\varepsilon_{\sigma}$ of $\sigma$ is 
defined in the usual manner via transpositions but any simple transposition 
$(i,i+1)$ of neighboring indices has the signature 
$(-1)^{(p_i-1)(p_{i+1}-1)}$.

There are two differentials on $B\A$: the inner differential
$$
d_{\A}([a_1|a_2|\dots|a_n])=\sum^n_{i=1}(-1)^i[\bar{a}_1|\bar{a}_2|\dots
|\bar{a}_{i-1}|d\,a_i|a_{i+1}|\dots|a_n],
$$
and the combinatorial differential  
$$
\delta([a_1|a_2|\dots|a_n])=
\sum^n_{i=1}(-1)^i[\bar{a}_1|\bar{a}_2|\dots|\bar{a}_{i-1}a_i
|a_{i+1}|\dots|a_n].
$$

These differentials are bihomogeneous and of bidegrees $(1,0)$ and $(0,1)$.
Moreover, 
$d_{\A}\delta+\delta\,d_{\A}=0$  and, therefore,
$\nabla=d_{\A}+\delta$ is a differential of bidegree $(1,1)$ on $B\A$. 
The Hopf algebra $B\A$ with the differential $\nabla$ is called
the {\it bar-construction} for $(\A,d)$.

The decreasing filtration
$\F_n(B\A)=\left(\sum_{j\geq n}\,B\A^{j,\ast}\right)$ is defined by the first
(tensor) grading and the spectral sequence corresponding to this filtration
is called the Eilenberg--Moore spectral sequence. It ``lives'' in the second 
quadrant. Although generally this filtration is infinite, the term
$\{E^{p,q}_{\infty}=\lim_{\stackrel{\to}{i}}\,E_i^{p,q}\}^{p,q}$    
is correctly defined and associated to $H(B(\A),\nabla)$. 
It is well-known \cite{M} that the bar-construction endowed only with
the combinatorial differential $\delta$ defines the homologies
of $\Tor^{\ast\ast}_{\A}(k,k)$ with coefficients in the $\A$-module $k$
where the $\A$-action is defined by the augmentation:
$a(k)=\varepsilon(a)\cdot k$.
It is easily derived that the term  $E_2$ of the Eilenberg--Moore
spectral sequence is
$$
E_2^{\ast,\ast}\simeq \Tor^{\ast}_{H^{\ast}(\A)}(k,k).
$$
The sequence is natural with respect to homomorphisms of differential gra\-ded
algebras. This means that for a given homomorphism 
$\varphi:(\A,d_{\A})\to(\B,d_{\B})$ the mapping of the bar-constructions
induces a homomorphism of the spectral sequences
$$
\varphi^{\ast}:(E_r(\A)^{\ast,\ast},d_r)\to(E_r^{\ast,\ast}(\B),d_r),
\ \ r=1,2,\dots.
$$
If $\A$ is of a topological origin, i. e., if
$(\A,d)=(\E^{\ast}(X),d)$ where $X$ is a smooth simply connected manifold, 
then for the corresponding Eilenberg--Moore spectral sequence we have

1) $E_2^{\ast,\ast}\simeq\Tor ^{\ast}_{H^{\ast}(X,\R)}(k,k)$;

2) $E_r\Rightarrow H^{\ast}(\Omega X,\R)$.

The cohomology suspension homomorphism is related to the Eilenberg--Moore
spectral sequence as follows. Consider a mapping $\Sigma:\A\to B(\A)$
defined as
\begin{equation}
\Sigma(a)=[a].
\label{1.6.1}
\end{equation}
This mapping commutes with the differentials and we obtain the induced
mapping
\begin{equation}
\Sigma^{\ast}:H^+(X)\to H(B(\A),\nabla)\simeq H^{\ast}(\Omega X,\R)
\label{1.6.2}
\end{equation}

It seems that the following proposition was first proved by Moore
\cite{JM} (see also \cite{LS}).

\begin{proposition}
\label{pr1.6.1}
There is an equality  
$$
\Sigma^{\ast}=\Omega^{\ast}.
$$
\end{proposition}

It follows from (\ref{1.6.1}) that $\Sigma^{\ast}$ descends through the 
Eilenberg--Moore spectral sequence and we have a sequence of mappings  
$$
\Sigma_r^{\ast}:H^+(X,\R)\to E_r^{-1,\ast}.
$$
In particular, $\Sigma_1^{\ast}:H^+(X,\R)\to E_1^{-1,\ast}=
B^{-1}(H^{\ast}(X,\R),\delta)$ is defined as 
$\Sigma_1^{\ast}(\alpha)=[\alpha]$,     
and there is the natural projection
$$
\Sigma_2^{\ast}:H^+(X,\R)\to E_2^{-1,\ast}=
\Tor^1_{H^{\ast}(X,\R)}(k,k)\simeq H^+(X,\R)/H^+\cdot H^+.
$$

\begin{corollary}
\label{cor1.6.2}
If $\rk h_{\ast} <\dim H^+(X,\R)/H^+\cdot H^+$, then
the corresponding Eilenberg--Moore spectral sequence  does not stabilize
at the term $E_2$.
\end{corollary}

\begin{center}
{\bf \S 3. Higher Massey operations in symplectic manifolds}
\end{center}

{\bf 3.1. Symplectic manifolds and symplectic blow-ups.}

A smooth manifold $X$ with a given $2$-form $\omega$ is called symplectic if
this form meets two conditions:

1) $\omega$ is nondegenerate, i. e., for any nonzero vector
$\xi\in T_x X$  there is a vector $\eta\in T_x X$ such that
$\omega(\xi,\eta)\ne 0$;

2) $d\,\omega=0.$

Such a form $\omega$ is called symplectic. 

By the Darboux theorem, near any point $x\in X$ there are local coordinates   
$x^1,\dots,x^{2n}$ such that the symplectic form equals      
$$
\omega=\sum_{j=1}^n d\,x^j\wedge d\,x^{j+n}.
$$
This implies that a symplectic manifold $X$ is even-dimensional and the 
form $\omega^n$ is the volume form:
$$
\omega^n=n!\,d\,x^1\wedge\dots\wedge d\,x^n.
$$
Take a Riemannian metric on $X$ and define an operator $A$ linearly 
acting in the fibers of the tangent bundle $TX$ as
$$
(A\xi,\eta)=g_{jk}(A\xi)^j\eta^k=\omega(\xi,\eta).
$$
The operator $A$ is skew-symmetric and, therefore, $A^{\ast}A=-A^2$ is 
symmetric and positive. Take the positive square root of $A^{\ast}A$:
$$
Q=\sqrt{-A^2}>0
$$
and put   $J=AQ^{-1}$. Then we have
$$
J^2=-1
$$
and  $J$ defines an almost complex structure on $X$. This structure is
compatible with the symplectic structure, which means that the form
$$
\langle \xi,\eta\rangle=(\xi,J\eta)
$$
is Hermitian and positive-definite.

In the sequel we shall consider only compact symplectic manifolds because
any noncompact almost complex manifold has a compatible symplectic
structure \cite{Gromov1}. For compact manifolds the condition  
$[\omega^n]\ne 0 \in H^{2n}(X,\R)$ shows that not any compact complex 
manifold has a symplectic structure.

Examples of symplectic manifolds are K\"ahler manifolds, .i e., complex 
manifolds with a Hermitian metric $h_{jk}\,d\,z^jd\,\bar{z}^k$ such that
the form
$$
\omega=h_{jk}\,d\,z^j \wedge d\,\bar{z}^k
$$
is symplectic.

Other examples of simply connected symplectic manifolds of dimension 
greater than four are constructed from known ones by the symplectic fibration
\cite{Th}, the symplectic blow-up \cite{Gromov1,McD}, and the fiber connected
sum \cite{G}. We shall consider symplectic blow-ups and first recall their 
construction.         

{\it Symplectic blow-up.}  Let $f:Y\to X$ be a symplectic embedding, i. e.,
$\omega$ be the symplectic form on $X$ and $f^{\ast}\omega$ be the symplectic
form on $Y$. Therefore we identify $Y$ with the submanifold  
$f(Y)\subset X$.

To every point $x\in Y\subset X$ there corresponds a subspace 
$E_x\subset T_x X$, which is the symplectic normal complement to $T_x Y$
in  $T_x X$:
$$
\omega(\xi,\eta)=0\ \ \ \ \xi\in T_x Y,\ \eta\in E_x.
$$
The spaces $E_x$ are pasted into a vector bundle   
$$
E\to Y,
$$
which is isomorphic to the normal bundle to $Y$.

Since $Y$ is a symplectic submanifold, the restrictions of
$\omega$ on $E_x$ are nondegenerate and there is a family of
almost complex structures $J_x$ on $E_x$ smoothly depending on $x$
and compatible with $\omega\big|_{E_x}$. Hence, the structure group of 
$E$ reduces to $U(k)=SO(2k)\cap Sp(k)$ where  
$k=\frac{1}{2}(\dim X-\dim Y)$.

We identify the fibers of  $E$ with  $\C^k$ and consider another fibration    
$$
\widetilde{E}\to Y,
$$
whose fibers are isomorphic to the fibers of the canonical line bundle      
$L\to \C P^{k-1}$ and which has the same associated principal      
$U(k)$-bundle as $E\to Y$.
The fibration
$$
L\to \C P^{k-1},
$$
where
$$
L=\{(\xi,l)\in\C^k\times\C P^{k-1}\big|\xi\in l\},
$$
parameterizes pairs $(l,\xi)$ where $l$ is a line in $\C^k$ and   
$\xi\in l$.    
The projection onto the base is
$$
(\xi,l)\to l\in\C P^{k-1}
$$
and the structure group $U(k)$ acts on the fibration as follows     
$$
A\cdot(\xi,l)=(A\xi,Al).
$$
The fibers of $E$ and, therefore, of $\widetilde{E}$
are endowed with the Hermitian metric  
$$
(\xi,\eta)=\omega(\xi,J\eta).
$$
We denote by $E_r$ and  $\widetilde{E}_r$ the submanifolds of 
$E$ and $\widetilde{E}$ defined by the inequality 
$|\xi|\le r$.

Let us construct now the symplectic blow-up $\widetilde{X}$ of $X$
along $Y$. For that take a closed tubular neighborhood $V$ of $Y$ and
identify it with $E_1$. Remove the interior of $V$ from $X$ and glue
$\widetilde{E}_1$ to the obtained manifold with boundary. 
This gluing is defined by the isomorphism
$\partial E_1=\partial\widetilde{E}_1=\partial V$.  
The resulted manifold
$$
\widetilde{X}=\overline{(X\backslash V)}\cup_{\partial V}\widetilde{E}_1
$$
is called the symplectic blow-up. In fact, it is a result of the fiberwise
blow-ups of the fibers $E_x$ at the points  $\xi=0$.

There is the projection 
$$
\pi:\widetilde{X}\to X,
$$
which is a diffeomorphism outside $Y$  
$$
\widetilde{X}\backslash\pi^{-1}(Y)\to X\backslash Y
$$
and its restriction over $Y$ is the fibration   
$$
\pi^{-1}(Y)\to Y
$$
with  $\C P^{k-1}$-fibers.

As shown in \cite{Gromov1} and \cite{McD},  there is a
symplectic structure on $\widetilde{X}$, which coincides with  
$\pi^{\ast}(\omega)$ outside a small neighborhood of     
$\pi^{-1}(Y)$.

{\bf 3.2. Survival of higher Massey products.}

The main aim of this subsection is a proof of the following statement.

\smallskip

\noindent

{\bf Theorem A} {\it Let a simply connected symplectic manifold $X$ have 
an irreducible generalized Massey product of dimension $k$. Then for any
symplectic submanifold $Y \subset X$ with
$\codim Y>k$ the corresponding symplectic 
blow-up $\widetilde{X}$ also has an irreducible generalized Massey product of
dimension $k$.}

\begin{corollary}
\label{cor2.2.1}
Under the hypothesis of the Theorem the manifolds $X$ and 
$\widetilde{X}$ are nonformal.
\end{corollary}

It follows from Theorem  A and Corollary \ref{cor1.5.4} that 

\begin{corollary}
\label{cor2.2.2}
Let the cohomology ring of a symplectic manifold
$X$ have a nonspheric generator of dimension $k$.
Then the symplectic blow-up $\widetilde{X}$ of $X$ along any symplectic 
submanifold $Y\subset X$, with $\codim Y >k$, also has a nonspheric 
$k$-dimensional generator and, therefore, is nonformal.
\end{corollary}

From that using Corollary \ref{cor1.5.2} we obtain

\begin{corollary}
\label{cor2.2.3}
Let a simply connected symplectic manifold $X$  
satisfies the inequality
$$
\rk h_k<\dim H^k(X)/H^k(X)\cap (H^+(X))^2.
$$
Then the symplectic blow-up $\widetilde{X}$ of $X$ along any symplectic
submanifold  $Y\subset X$, with $\codim Y > k$, is nonformal.
\end{corollary}

Let us prove Theorem A.
Let $\pi:\widetilde{X}\to X$ be the natural projection.
By assumption, there is a generalized irreducible Massey product
$u\in H^k(X)$. We shall show that   
$\pi^{\ast}(u)\in H^k(\widetilde{X})$ 
is also  a generalized irreducible Massey product.

Let $A\in M(\E^{\ast}(X))$ be a solution to the Maurer--Cartan equation such 
that $u$ is an entry of the matrix $[d\,A-\bar{A}\cdot A]$.
Since solutions to the Maurer--Cartan equation are natural,
$\pi^{\ast}(A)$ is a solution to this equation in
$\E^{\ast}(\widetilde{X})$. Moreover the corresponding entry of the matrix
$$
[d\,\pi^{\ast}(A)-\overline{\pi^{\ast}(A)}\cdot \pi^{\ast}(A)]=
\pi^{\ast}([d\,A-\bar{A}\cdot A])
$$
equals $\pi^{\ast}(u)$. It remains to prove that $\pi^{\ast}(u)$ is
irreducible. This is implied by the following statement.

\begin{lemma}
\label{lem2.2.4}
Let $\pi:\widetilde{X}\to X$ be the natural projection of the
symplectic blow-up of $X$ along $Y\subset X$. The induced
mapping  
\begin{equation}
\pi^{\ast}:\Tor^{H^{\ast}(X)}_{1,k}(\R,\R)\to
\Tor^{H^{\ast}(\widetilde{X})}_{1,k}(\R,\R)
\label{2.2.1}
\end{equation}
is a monomorphism for $k<\codim Y$.
\end{lemma}

Indeed, for the algebra $H=H^{\ast}(X)$ there is an isomorphism   
$$
\Tor^H_{1,\ast}(k,k)\simeq H^+/H^+\cdot H^+
$$
and the monomorphism (\ref{2.2.1}) implies that 
$\pi^{\ast}(u)\neq 0\, \mod\, H^+(\widetilde{X})\cdot H^+(\widetilde{X})$,
which means that $\pi^{\ast}(u)$ is irreducible.

{\sl Proof of the lemma.}
Let $V\subset X$ be a tubular neighborhood of $Y$ and let
$\widetilde{V}$ be its symplectic blow-up along $Y$.
In this notation,
we have
$\widetilde{X}=(X\backslash\dot{V})\cup_{\partial\widetilde{V}}\,
\widetilde{V}$ where $\dot{V}=V\backslash\partial V$ and the boundaries of
$\partial V$ and  $\partial\widetilde{V}$ are identified in the natural 
manner. Consider the space 
$$
\widehat{X}=X\cup_{\partial\widetilde{V}}\,
\widetilde{V}=\widetilde{X}\cup_{\partial V}\, V.
$$
There is a commutative diagram   
\begin{equation}
\label{2.2.2}
\includegraphics{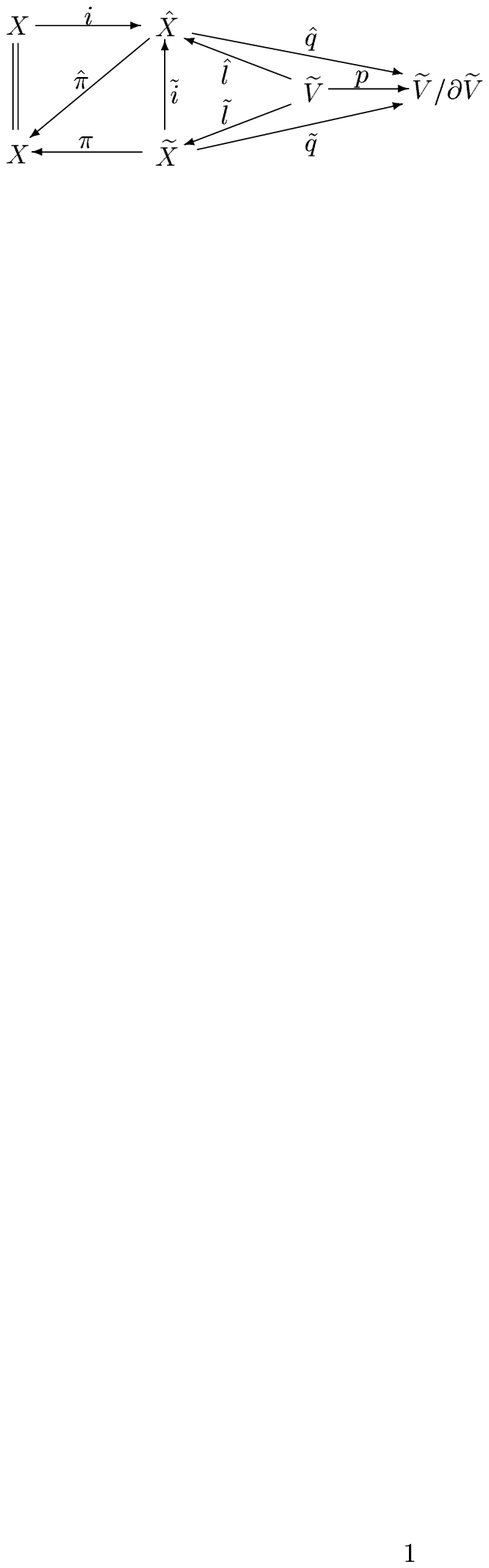}
\end{equation}
where $i$, $\widetilde{i}$, $\widehat{l}$, and $\widetilde{l}$ are the
embeddings, $p$, $\widehat{q}$, $\widetilde{q}$, $\pi$, and 
$\widehat{\pi}$ are the projections and, moreover, 
$\widehat{\pi}$ is a retraction.

Consider the ideal  $I=\Ker i^{\ast}$. Since $X$ is a retract of  
$\widehat{X}$, there is a multiplicative isomorphism  
\begin{equation}
H^{\ast}(\widehat{X})\simeq H^{\ast}(X)\oplus I
\label{2.2.3}
\end{equation}
where $H^{\ast}(X)$ is identified with its image under the 
retraction $\widehat{\pi}^{\ast}$. 
Since $I$ is an ideal, by (\ref{2.2.3}), we conclude that 
\begin{equation}
\widehat{\pi}^{\ast}:H^+(X)/(H^+(X))^2\longrightarrow 
H^+(\widehat{X})/(H^+(\widehat{X}))^2
\label{2.2.4}
\end{equation}
is a monomorphism.

Now from the exact cohomology 
sequence for $(\widehat{X},X)$ derive that
$\widehat{q}^{\ast}$ maps $H^+(\widetilde{V}/\partial\widetilde{V})$
isomorphically onto $I$. Let 
$\widetilde{a}\in H^2(\widetilde{V}/\partial\widetilde{V})$ be the projective 
Thom class. Using the fibration
\begin{equation}
r:\widetilde{V}\stackrel{\dot{\C P}^m}{\longrightarrow} Y
\label{2.2.5}
\end{equation}
where $r$ is the projection, $m=\frac{1}{2}\,\codim Y$, and  
$\dot{\C P}^m=\C P^m\backslash D^{2m}$, we describe  
$H^+(\widetilde{V}/\partial\widetilde{V})$ as a $H^{\ast}(Y)$-module:
$$
H^+(\widetilde{V}/\partial\widetilde{V})\simeq\widetilde{a} 
H^{\ast}(Y)\oplus\widetilde{a}^2 H^{\ast}(Y)
\oplus\dots\oplus\widetilde{a}^m H^{\ast}(Y).
$$
Put $\widehat{a}=\widehat{q}^{\ast}(\widetilde{a})$ and applying 
$\widehat{q}^{\ast}$ we obtain
\begin{equation}
I\simeq\widehat{a} H^{\ast}(Y)\oplus\widehat{a}^2 
H^{\ast}(Y)\oplus\dots\oplus\widehat{a}^m H^{\ast}(Y).
\label{2.2.6}
\end{equation}
The fibration (\ref{2.2.5}) defines a structure of a $H^{\ast}(Y)$-module 
on $H^{\ast}(\widetilde{V})$:  
\begin{equation}
H^{\ast}(\widetilde{V})\simeq H^{\ast}(Y)\oplus 
a H^{\ast}(Y)\oplus\dots\oplus a^{m-1}H^{\ast}(Y)
\label{2.2.7}
\end{equation}
where  $a\in H^2(\widetilde{V})$ is the two-class, which is mapped into the
generator $\dot{\C P}^m\subset\widetilde{V}$. Notice that
$a=\widehat{l}^{\ast}(\widehat{a})=p^{\ast}(\widetilde{a})$ and 
$\widetilde{a}=\widehat{q}^{\ast}(\widehat{a})$. 
Moreover the decompositions (\ref{2.2.6}) and (\ref{2.2.7}) are 
related by the homomorphism induced by the embedding  
$\widehat{l}$. Finally we put
$$
I_1=\widehat{a} H^{\ast}(Y)\oplus\dots\oplus\widehat{a}^{m-1} H^{\ast}(Y),
\ \ I_2=\widehat{a}^m H^{\ast}(Y)
$$
and show that $\widetilde{i}^{\ast}\big |_{I_1}$ is a monomorphism.
Define a $H^{\ast}(Y)$-submodule $H^{\ast}(\widetilde{V})$ as follows 
$$
H_+^{\ast}(\widetilde{V})=a H^{\ast}(Y)\oplus a^2 H^{\ast}(Y)\oplus
\dots\oplus a^{m-1} H^{\ast}(Y).
$$
As shown in \cite{McD} there is an exact short sequence of modules
\begin{equation}
0\to H^{\ast}(X)\stackrel{\pi^{\ast}}{\to} H^{\ast}(\widetilde{X})
\stackrel{\widetilde{l}^{\ast}}{\to} H_+^{\ast}(\widetilde{V})\to 0.
\label{2.2.8}
\end{equation}
It follows from the diagram (\ref{2.2.2}) that 
$\widehat{l}^{\ast}\big |_{I_1}$ is an isomorphism and
$\widetilde{i}^{\ast}\big |_{I_1}$ is a monomorphism, which splits
$\widetilde{l}^{\ast}$.

Now using the decomposition (\ref{2.2.3}) and keeping in mind that
$\widetilde{i}^{\ast}$ is a multiplicative isomorphism in dimensions less or 
equal than $\codim Y$ we obtain 
$$
(H^+(\widetilde{X})\cdot H^+(\widetilde{X}))^k\subseteq(H^+(X)\cdot H^+(X))
\oplus(i^{\ast}(I))^k,\ \
k<\codim Y.
$$
From that we conclude that 
$$
\pi^{\ast}:H^k(X)/H^k(X)\cap(H^+(X))^2\to 
H^k(\widetilde{X})/H^k(\widetilde{X})
\cap (H^+(\widetilde{X}))^2
$$
is a monomorphism for  $k<\codim Y$.  This proves the lemma.

{\sl Remark.} For the simplest case, for triple Massey products, the notions of
nontriviality and irreducibility are very closed. By a simple modification
of the given proof of Lemma \ref{lem2.2.4} one can obtain the following
statement.     

\begin{proposition}
\label{pr2.2.5}
Let  $X$ be a simply connected symplectic manifold with a nontrivial triple
Massey product of dimension $k$. Then for any its symplectic blow-up 
$\widetilde{X}$ along
$Y\subset X$, with $\codim Y > k$, also has a nontrivial Massey product of 
dimension $3$.
\end{proposition}

{\bf 3.3. Symplectic manifolds of arbitrary large weight}

Let  $\langle X_1,\dots,X_k\rangle$ be a $k$-tuple matrix Massey product
in $H^{\ast}(M)$. We say that the weight of
$\langle X_1,\dots,X_k\rangle$ is {\it strictly equal to $k$}
if there is an element $a\in\langle X_1,\dots,X_k\rangle$, which does not 
belong to the linear subspace generated by all Massey products 
$\langle Y_1,\dots,Y_l\rangle$ with $l<k$. We shall see below (Remark 3) that
nontrivial relations between cohomologies classes sometimes may 
imply some nontrivial $k$-tuple products are expressed in terms of
Massey products of less weight and of another elements. 
The main result of this subsection is the following

\begin{theorem}
\label{th4.2.1}
For any $k$ there exist symplectic manifolds with Massey products whose 
weights are strictly equal to $2k$.
\end{theorem}

For proving this theorem we shall use the symplectic nilmanifolds $M(2n)$
introduced by us in \cite{BT}. Let us recall in brief their construction.

Take the Witt algebra $W(1)$ of formal vector fields on the line.
It has a basis  
$$
e_i=x^{i+1}\frac{d}{dx},\ \ i=-1,0,1,\dots,
$$
with the commutation relations     
\begin{equation}
[e_i,e_j]=(j-i)e_{i+j},\  \ i,j\geq -1.
\label{4.2.1}
\end{equation}
Denoting by $\L(\dots)$ the linear span of the corresponding vectors, we
consider the subalgebras
$$
\L_k(1)=\L(e_k,e_{k+1},\dots),\ \ k=1,2,\dots.
$$
There is a natural filtration
$$
\dots\subset\L_1(1)\subset\L_0(1)\subset\L_{-1}=W(1),
$$
which is useful for study $W(1)$ and its subalgebras and quotients.

Consider the finite-dimensional nilpotent Lie algebras
$$
\V_n=\L_1(1)/\L_{n+1}(1), \ \ n=1,2,\dots,
$$
and denote the corresponding nilpotent Lie groups by $V_n$, $n=1,2,\dots$.  
By (\ref{4.2.1}), $\V_n$ is a nilpotent $n$-dimensional Lie algebra 
with a basis $\{e_1,\dots,e_n\}$ and the Lie brackets
\begin{equation}
[e_i,e_j]=\cases{
            (j-i)e_{i+j}, & $i+j\leq n$,\cr
            0, &  $i+j>n$.
           }
\label{4.2.2}
\end{equation}
By the nilpotence, we identify $\V_n$ and $V_n$ as sets and
define the multiplication $\times$ by the Campbell--Hausdorff formula.
Another interpretation of $V_n$ as a group of polynomial transformations of 
the line was recently used in \cite{Buh}.

The formulas (\ref{4.2.2}) show that the structure constants of $\V_n$ 
are rational and, therefore, by the Mal'tsev theorem \cite{Mal}, 
$V_n$ has uniform lattices. Among such lattices we choose the group 
$\Gamma_n$ generated by $\{e_1,\dots,e_n\}$ with respect to $\times$.
We obtain a family of nilmanifolds
$$
M(n)=V_n/\Gamma_n,\ \ n=1,2,\dots.
$$
Notice that $V_3$ is the three-dimensional Heisenberg group $\H$, 
$M(3)\simeq\H/\H_{\Z}$, and 
$M(3)\times S^1$ is the Kodaira--Thurston manifold.

The nilmanifolds $M(n)$ with $n\geq 3$ have many prominent properties and,
in particular, for even $n$ they have natural symplectic structures \cite{BT}. 
Let $\{\omega_1,\dots,\omega_n\}$ be the basis for left-invariant $1$-forms  
on$V_n$ dual to $\{e_1,\dots,e_n\}$.  It follows from (\ref{4.2.2}) that 
\begin{equation}
d\,\omega_k=(k-2)\,\omega_1\wedge\omega_{k-1}+(k-4)\,
\omega_2\wedge\omega_{k-2}+\dots.
\label{4.2.3}
\end{equation}
Moreover, as shown in \cite{BT} the left-invariant $2$-form
$$
\Omega_{2m}=(2m-1)\,\omega_1\wedge\omega_{2m}+(2m-3)\,
\omega_2\wedge\omega_{2m-1}+\dots+
\omega_m\wedge\omega_{m+1}
$$
defines a left-invariant symplectic structure on $V_{2m}$.
The factorization descends this structure to $M(2m)$ where it defines
an integer symplectic form, which we denote also by $\Omega_{2m}$.  
Let $[\Omega_{2m}]\in H^2(M(2m),\Z)$ the cohomology class of this form.
Theorem \ref{th4.2.1} is implied by the following statement.

\begin{proposition}
\label{pr4.2.2}
The symplectic class  $[\Omega_{2m}]\in H^2(M(2m),\Z)$ is a $2m$-tuple
matrix Massey product whose weight is strictly equal $2m$. 
\end{proposition}

{\sl Proof.}  
Consider an exterior algebra 
$\Lambda_{2m}=\Lambda(\omega_1,\omega_2,\dots,\omega_{2m};d)$ where 
$d\omega_1=d\omega_2=0$ and 
$d\omega_k$ are defined by the formula (\ref{4.2.3}) for $k\geq 3$. 
By the Nomizu theorem \cite{N}, 
the embedding
$$
\Lambda_{2m}\longrightarrow\E^{\ast}(M(2m)),
$$
induces by the embedding of left-invariant forms into the complex of all forms
on $M(2m)$ is a weak isomorphism, Hence, we may compute Massey products
in $H^{\ast}(M(2m))$ from $\Lambda_{2m}$. Notice that 
$\Lambda_{2m}$ is isomorphic to the minimal model for $M(2m)$.

Let us represent $\Omega_{2m}$ by a $2m$-tuple Massey product.  
There is an inclusion
\begin{equation}
\begin{array}{c}
[\Omega_{2m}]\in\langle-m\omega_1,-(m-1)\omega_1,\dots,
-2\omega_1,(-\omega_1\;0),
\left(\begin{array}{c} \omega_2 \\ \omega_1 \end{array}\right), \\
2\omega_1,\dots,
(m-1)\omega_1,\omega_1\rangle.
\end{array}
\label{4.2.4}
\end{equation}

For proving (\ref{4.2.4}) consider the 
$(2m+2)\times(2m+2)$-matrix $A$ of the form 
$$
\left(\begin{array}{cccccccc}
0 & -m\omega_1 & (1-m)\omega_2 &\dots & (k-m)\omega_{1+k} & \dots & 
(m-1)\omega_{2m} & 0 \\
0 & 0          & (1-m)\omega_1 &\dots & (k-m)\omega_k & \dots & 
(m-1)\omega_{2m-1} & \omega_{2m} \\
& \dots & & & & & \dots & \\
0 & 0 & 0 & \dots & 0 & \dots & (m-1)\omega_1 & \omega_2  \\
0 & 0 & 0 & \dots & 0 & \dots & 0 & \omega_1 \\
0 & 0 & 0 & \dots & 0 & \dots & 0 & 0
\end{array}\right).
$$
Straightforward computations show that $A$ satisfies the Maurer--Cartan
eq\-uation and, therefore, is a defining system for the product 
(\ref{4.2.4}). These computations also show that 
the matrix $\bar{A}\wedge A-d\,A$ has one and only one nonzero entry.
It is the $(1,2m+2)$ entry, which equals
$$
-m\bar{\omega}_1\wedge\omega_{2m}-(m-1)\bar{\omega}_2\wedge
\omega_{2m-1}-\dots-
\bar{\omega_m}\wedge\omega_{m+1}+\bar{\omega}_{m+2}\wedge\omega_{m-1}+\dots+
$$
$$
+(m-1)\bar{\omega}_{2m}\wedge\omega_1=
\sum_{i=1}^m[2(m-i)+1]\omega_i\wedge\omega_{2m-i+1}
=\Omega_{2m}.
$$
The last equality proves the inclusion (\ref{4.2.4}).

Now show that $\Omega_{2m}$ can not be expressed as a linear combination of
$k$-tuple  Massey products with $k<2m$. For that we notice that the algebra
$\Lambda_{2m}$ is bigraded. The gradings are defined on the generators
$$
\deg_1(\omega_i)=1,\ \ \deg_2(\omega_i)=i,\ \ i=1,2,\dots,2m,
$$
and are extended by the multiplication. It follows from  (\ref{4.2.3}) that
the bidegree of the differential $d$ is $(1,0)$. Notice, that the degree rules
used for the definition of multituple Massey products and introduced in
\S 2.3 are applied only to $\deg_1$.

The algebra $\Lambda_{2m}$ also has an increasing filtration    
$$
0=\F_0\subset\F_1\subset\dots\subset\F_{2m}=\Lambda_{2m}
$$
where $\F_k=\Lambda(\omega_1,\dots,\omega_k)$, $k=1,2,\dots,2m$.   
This filtration meets the following conditions:

(1) $f(\alpha\wedge\beta)\leq\max(f(\alpha),f(\beta))$.

(2) $f(d\,\alpha)<f(\alpha)$, and, moreover, if $\deg_1(\alpha)=1$, then
$f(d\,\alpha)=f(\alpha)-1$. \\
The last property follows from (\ref{4.2.3}).

Suppose that
\begin{equation}
\Omega_{2m}=\sum_s\alpha_s\langle X_1^{(s)},
\dots,X_{k(s)}^{(s)}\rangle,\ \ k(s)<2m,\ \
\alpha_s\in\R
\label{4.2.6}
\end{equation}
and  $X_j^{(s)}$ are matrices whose entries are 
closed elements from $\Lambda_{2m}$. Since 
$\deg_1\Omega=2$, the degree rules imply that $\deg_1 X_j^{(s)}=1$. 
Since $X_j^{(s)}$ is closed, all its entries are linear combinations of
$\omega_1$  and $\omega_2$ and, therefore,
\begin{equation}
f(X_j^{(s)})\leq 2\ \ \mbox{for all}\ \ s\ \mbox{and}\ \ j.
\label{4.2.7}
\end{equation}

Let $A^{(s)}=(X_{(i,j)}^{(s)})$ be a defining system for
$\langle X_1^{(s)},\dots,X_{k(s)}^{(s)}\rangle$.    
Applying the degree rules, we obtain 
$\deg_1 X_{(i,j)}^{(s)}=1$ for all  $s$, $i$, $j$. 
Hence, it follows from (\ref{4.2.7}) and  (2)  that 
$f(X_{(i,i+1)}^{(s)})\leq 3$.  Now applying the equality 3 from Definition
\ref{def1.3.4} we obtain, by induction, that   
$$
f(X_{(i,j)}^{(s)})\leq j-i+2.
$$
The last step of this inductive process together with the formula 4 from
Definition \ref{def1.3.4} implies
$$
f(c(A^{(s)}))\leq k(s)
$$
where $c(A^{(s)})$ is the cocycle of the defining system.
But from another side we have
$f(\Omega(2m))$ $=2m$ and, therefore, the equality
(\ref{4.2.6}) does not hold if $k(s)<2m$ for all $s$.
This contradiction proves the proposition.

{\sl Remarks.} 1. We prove Proposition \ref{pr4.2.2} by presenting 
an explicit formula (\ref{4.2.4}) for the product representing
$[\Omega_{2m}]$. This product is not uniquely defined. For instance,
for $M(4)$ there is another inclusion except (\ref{4.2.4}):
$$
[\Omega_4]\in\langle(0\ \ 6\omega_2),
\left(\begin{array}{cc}\omega_1 & \omega_2 \\ 0 & \omega_1\end{array}\right),
\left(\begin{array}{cc}\omega_1 & \omega_2 \\ 0 & \omega_1\end{array}\right),
\left(\begin{array}{c}\omega_2 \\ \omega_1\end{array}\right)\rangle.
$$
We are leave that for the reader to find the corresponding solution to the
Maurer--Cartan equation.

2. In  $M(4)$ the symplectic class is a sum of two cohomology classes        
$$
[\Omega_4]=[3\omega_1\wedge\omega_4]+[\omega_2\wedge\omega_3].
$$
It is easy to check that the first component is represented by
a quadruple Massey product as follows
\begin{equation}
[3\omega_1\wedge\omega_4]\in\langle 6\omega_2,
\omega_1,\omega_1,\omega_1\rangle\label{4.2.8}
\end{equation}
and its weight is strictly equal to $4$. Another component is 
represented in the form       
\begin{equation}
[\omega_2\wedge\omega_3]\in\langle -\omega_1,\omega_2,\omega_2\rangle
\label{4.2.9}
\end{equation}
and the weight of this product is strictly equal to $3$.

3. For $m>2$ the both classes  
$[3\omega_1\wedge\omega_4]$  and $[\omega_2\wedge\omega_3]$ 
do not vanish in  $H^2(M(2m))$, but, by the equality    
$$
d\,\omega_5=\Omega_4=3\omega_1\wedge\omega_4+\omega_2\wedge\omega_3,
$$
their sum equals zero. Hence, for $m >2$ there is a nontrivial relation
between the quadruple product (\ref{4.2.8}) and the triple product
(\ref{4.2.9}) in  $H^{\ast}(M(2m))$. Therefore, for $m >2$
$[3\omega_1\wedge\omega_4]\in H^2(M(2m))$ is not a cohomology class
whose weight is strictly equal to $4$.

\smallskip

{\bf 3.4. Inheritance of higher Massey products.}

Let $(X,\omega)$ be a symplectic manifold and let $Y\subset X$ be its 
submanifold. 

In this subsection we find sufficient conditions, which guarantee that
some nontrivial higher matrix Massey products in $H^{\ast}(Y)$ are 
inherited after the 
blow-up, i.e., generate nontrivial matrix Massey products in       
$H^{\ast}(\widetilde{X})$. 

Consider matrix cohomology classes  
$S_1,\dots,S_n\in N(H^+(\widetilde{X}))$ and assume that the product
$\langle S_1,\dots,S_n\rangle$ is defined and nontrivial. Our aim 
is to find some corresponding $n$-tuple Massey product 
$\langle\widetilde{S}_1,\dots,\widetilde{S}_n\rangle$, with
$\widetilde{S}_i\in N(H^+(Y))$,  which is nontrivial in 
the symplectic blow-up  of $X$ along $Y\subset X$. 
If $n \geq 5$ and matrix classes $S_i$, $1\leq i\leq n$, are taken arbitrary,
then it is already a problematic question can one effectively achieve that 
for $Y\subset \C P^N$ with arbitrary large $N$.
The reasons for that are the structure of Massey products for $n \geq 5$ and
possible relations between Massey products of different orders in
$H^{\ast}(Y)$. But there are two important cases, which we shall consider in 
the sequel:

1)  $S_i$ are matrices of one-dimensional classes, i.e., 
$S_i\in N(H^1(Y))$, $1\leq i\leq n$, where $n$ is arbitrary;

2)  $S_i$ are arbitrary matrices for $1\leq i\leq n$ and $n=3,4$.

Let us formulate the main results of this subsection.

\smallskip

\noindent
{\bf Theorem B} {\it Let a symplectic manifold $(Y,\omega)$ have a 
nontrivial 
matrix $n$-tuple Massey product $\langle S_1,\dots,S_n\rangle$ where
$S_i\in N(H^1(Y))$ are matrices of one-dimensional
cohomology classes for $1\leq i\leq n$. Then for 
any symplectic embedding $Y\subset X$ of codimension
not less than $2(n+1)$ the corresponding symplectic blow-up 
$\widetilde{X}$ has a nontrivial $n$-tuple Massey product
$\langle\widetilde{S}_1,\dots,\widetilde{S}_n\rangle$, where
$\widetilde{S}_i\in N(H^3(\widetilde{X}))$, $1\leq i\leq n$.}

\begin{corollary}
\label{cor2.4.1}
For any $k$ there exists a simply connected symplectic manifold $X$ of
dimension $6k+2$ with a nontrivial $2k$-tuple matrix Massey product in
$H^{4k+2}(X)$.
\end{corollary}

Indeed, consider the symplectic manifold $M(2k)$ defined in \S 3.3.
Let $f:M(2k)\to\C P^{3k+1}$  be a symplectic embedding, which exists by the 
Gromov--Tischler theorem \cite{T}. The corollary follows from Proposition
\ref{pr4.2.2} and Theorem  B applied to the pair   
$M(2k)\subset\C P^{3k+1}$. Corollary \ref{cor2.4.1} is established.

\smallskip

\noindent
{\bf Theorem C} {\it Let a symplectic manifold $(Y,\omega)$ have a 
nontrivial triple matrix Massey product. Then for any symplectic 
embedding $Y\subset X$
of codimension greater or equal than $8$ the corresponding symplectic blow
up $\widetilde{X}$ also has a nontrivial triple matrix Massey product.}

In the sequel the following definition will be useful.    

\begin{definition}
\label{def2.4.2}
Let $\A$ be a differential graded algebra and let $H\in N(\A)$. The upper 
degree of $H$ is
$$
\sdeg H=\max_{i,j}\deg h_{ij}
$$
where  $h_{ij}\in\A$ are the matrix entries of $H$.
\end{definition}

\smallskip

\noindent
{\bf Theorem D} {\it Let  a symplectic manifold $(Y,\omega)$ have a 
strictly irreducible quadruple matrix Massey product 
$\langle S_1, S_2,S_3,S_4\rangle$. Then for any symplectic submanifold 
$Y\subset X$ such that
$$
\codim Y>2\,\sdeg\langle S_1,S_2,S_3,S_4\rangle
$$
the corresponding symplectic blow-up $\widetilde{X}$ has a nontrivial 
quadruple matrix Massey product.}

Let us prove these theorems. First  we 
construct from the given Massey product   
$\langle S_1,\dots,S_n\rangle$ in $H^{\ast}(Y)$
another Massey product 
$\langle\widetilde{S}_1,\dots,\widetilde{S}_n\rangle$ in
$H^{\ast}(\widetilde{X})$. These parts of the proofs are similar for all 
theorems.

Let $S_1,\dots,S_n\in N(H^{\ast}(Y))$ be multipliable matrices and 
the Massey product $\langle S_1,\dots,S_n\rangle$ be defined.  
If $\pi:\widetilde{V}\to V\sim Y$ is the natural projection, then  
$\pi^{\ast}:H^{\ast}(Y)\to H^{\ast}(\widetilde{V})$ is a monomorphism and,
hence,  
we shall not differ $S_i$ and $\pi^{\ast}(S_i)$ for $1\leq i\leq n$.
Consider the matrix classes $a\cup S_i\in N(H^{\ast}(\widetilde{V}))$, 
$1\leq i\leq n$, where $a\in H^2(\widetilde{V})$ is the same as in \S 3.2.
It is shown in \S 2.2 that
$$
p^{\ast}:H^{\ast}(\widetilde{V}/\partial\widetilde{V})\to 
H^{\ast}(\widetilde{V})
$$
maps $\widetilde{a}\otimes H^{\ast}(Y)
\subset H^{\ast}(\widetilde{V}/\partial\widetilde{V})$ isomorphically
onto the corresponding submodule    
$a\cup H^{\ast}(Y)\subset H^{\ast}(\widetilde{V})$
where $p^{\ast}(\widetilde{a})=a$.  
Therefore, the elements  
$S^{\prime}_i\in N(H^{\ast}(\widetilde{V}))$ 
such that $p^{\ast}(S^{\prime}_i)=a\cup S_i$ 
are correctly defined for $1\leq i\leq n$.
Finally, we put
\begin{equation}
\widetilde{S}_i=\widetilde{q}^{\ast}(S^{\prime}_i),\ \ 1\leq i\leq n
\label{2.4.1}
\end{equation}
and notice that, since the diagram on (\ref{2.2.2}) commutes, we have
$$
\widetilde{l}^{\ast}(\widetilde{S}_i)=a\cup S_i=p^{\ast}(S^{\prime}_i),
\ \ 1\leq i\leq n.
$$
The theorems are implied by the following statement.

\begin{lemma}
\label{lem2.4.3}
Under the hypotheses of Theorems B, C, and D the matrix Massey 
$\langle\widetilde{S}_1,\dots,\widetilde{S}_n\rangle$ is defined and
nontrivial.
\end{lemma}

Since Massey products are natural and 
the diagram (\ref{2.2.2}) is commutative, Lemma \ref{lem2.4.3}
follows from these two lemmas. 

\begin{lemma}
\label{lem2.4.4}
Under the hypotheses of Theorems B, C, and D the Massey product
$\langle S^{\prime}_1,\dots,S^{\prime}_n\rangle$ is defined.
\end{lemma}

\begin{lemma}
\label{lem2.4.5}
Under the hypotheses of Theorems B, C, and D the Massey product
$\langle a\cup S_1,\dots,a\cup S_n\rangle$ is nontrivial.
\end{lemma}

Notice that it is clear that the latter product is defined.

{\sl Proof of Lemma  \ref{lem2.4.4}.} 
Let  $A=(X(i,j))$, where $1\leq i\leq j\leq n$ and$(i,j)\neq (1,n)$, be 
a defining system for $\langle S_1,\dots,S_n\rangle$. 
We have $X(i,j)\in N(\E^+(Y))$ and $[X(i,i)]=S_i$ for $1\leq i\leq n$.   
Denote by $X^{\prime}(i,j)$ the matrix differential forms 
$\pi^{\ast}(X(i,j))$ where $1\leq i\leq j\leq n$.
Let $\alpha$ be a  $2$-form representing the class $a\in H^2(\widetilde{V})$.
It can be taken such that $\alpha\big|_U=0$ for $U$, 
a small tubular neighborhood of $\partial\widetilde{V}$.
Indeed, $j^{\ast}(a)=0$ where $j:\partial\widetilde{V}\to\widetilde{V}$ 
is the embedding and this means that $\alpha=d\,\beta$ in $U^{\prime}$
for some tubular neighborhood $U^{\prime}\supset U$, which contains $U$. 
Take a cut-off function $\varphi$, which equals $0$ in  $U$ and equals $1$
outside $U^{\prime}$ and consider a form $\alpha_1$ such that 
$\alpha_1=d\,(\varphi\beta)$ in $U^{\prime}$ and $\alpha_1=\alpha$
outside $U^{\prime}$. Since $\alpha_1\sim \alpha$, we obtain a form 
representing $a$ and with the necessary properties.

It is clear that the family of matrix forms
$$
(X^{\prime\prime}(i,j)=\alpha^{j-i+1}\wedge X^{\prime}(i,j)),
\ \ 1\leq i\leq j\leq n,
$$
form a defining system for     
$\langle a\cup S_1,\dots,a\cup S_n\rangle$. 
Since $X^{\prime\prime}(i,j)\big|_U\equiv 0$ for $1\le i\le j\le n$
these forms may be extended by zero forms onto      
$\widetilde{V}\cup_{\partial\widetilde{V}} C\partial\widetilde{V}\simeq
\widetilde{V}/\partial\widetilde{V}$.        
Hence we obtain a family of simplicial forms on 
$\widetilde{V}/\partial\widetilde{V}$. We preserve for these 
forms the notation
$A^{\prime\prime}=(X^{\prime\prime}(i,j))$, 
$1\le i\le j\le n$. By the construction, $A^{\prime\prime}$
is a defining system for     
$\langle S_1^{\prime},\dots,S_n^{\prime}\rangle$.   
This proves Lemma \ref{lem2.4.4}.

The proof of Lemma  \ref{lem2.4.5} splits into three cases, which imply
Theorems B, C, and  D.

{\bf Lemma \ref{lem2.4.5} B}
\label{lem2.4.5B}
{\it Under the hypothesis of Theorem B the matrix Massey product
$\langle a\cup S_1,\dots,a\cup S_n\rangle$ is nontrivial.}

{\sl Proof.}     
In difference with the previous proof it is more convenient to use in this 
case minimal algebras but not the algebras of forms.
Let $(\M_Y,d_Y)$ be the minimal model for $Y$ and let 
$m=\frac{1}{2}\,\codim Y$.    
Then the minimal model for $\widetilde{V}$ is as follows \cite{BT}:
$$
(\M,d)\simeq\M_Y\otimes_d\{x,y\}
$$
where  $\deg x=2$, $\deg y=2m-1$, $d\big|_{\M_Y}=d_Y$, $d\,x=0$,
\begin{equation}
d\,y=x^m+c_1\wedge x^{m-1}+c_2\wedge x^{m-2}+\dots+c_{m-1}\wedge x+c_m,
\label{2.4.2}
\end{equation}
$c_i\in\M_Y$, $1\le i\le m$. Notice that  $[x]=a$.

Choose a representative
$R_i\in N(\M_Y^+)$, $[R_i]=S_i$, for 
any cohomology class $S_i$ where $1\le i\le n$.
We have
$$
[x\wedge R_i]=a\cup S_i,\ \ 1\le i\le n.
$$
Suppose that the Massey product is trivial:
$$
0\in\langle a\cup S_1,\dots,a\cup S_n\rangle
$$
and  $A=(X(i,j))$, $1\le i\le j\le n$, is a defining system for the product.
Since the entries of the matrices $S_i$ are one-dimensional 
cohomology classes, we have
$$
\deg X(i,j)=2(j-i)+3,\ \ 1\le i\le j\le n,\ \ (i,j)\ne (1,n).
$$
Here each of the matrices contains entries of the same degree and hence
the degrees of matrices are correctly defined.

Expand  $X(i,j)$ into powers of $x$:
\begin{equation}
X(i,j)=\sum_{l=0}^{j-i}x^l\wedge X_l(i,j)+x^{j-i+1}\wedge X_{j-i+1}(i,j)
\label{2.4.3}
\end{equation}
where $X_l(i,j)\in\M_Y$. Since $X(i,i)=x\wedge R_i$ for $1\le i\le n$,
we have 
\begin{equation}
X_1(i,i)=R_i\ \ (X_0(i,i)=0), 1\le i\le n.\label{2.4.4}
\end{equation}
By the hypothesis of Theorem B, we have $j-i+1\le n-1<m$ for all $X(i,j)$.  
Therefore the equalities
$$
d\,X(i,j)=\sum_{r=i}^{j-1}\overline{X(i,r)}\wedge X(r+1,j)
$$
and (\ref{2.4.2}) imply that 
$$
d\,X_{j-i+1}(i,j)=\sum_{r=i}^{j-1}\overline{X_{r-i+1}(i,r)}\wedge 
X_{r+1}(r+1,j).
$$
The latter equality together with (\ref{2.4.4}) means that 
$A_1=(X_{j-i+1}(i,j))$, with $1\le i\le j\le n$, is a defining
system for $\langle S_1,\dots,S_n\rangle$.
Applying to $c(A)$ the analogous computation with taking 
(\ref{2.4.3}) into account, we obtain
$$
c(A)=\sum_{r=1}^{n-1}\overline{X(1,r)}\wedge X(r+1,n)=
$$
$$
=x^n\wedge\sum_{r=1}^{n-1}\overline{X_r(1,r)}\wedge 
X_{n-r}(r+1,n)+\dots=x^n\wedge c(A_1)+\dots,
$$
where dots denote terms whose degree in $x$ is less than $n$.

Since  $n<m$, the latter equality means that
$[c(A)]=0$ implies $[c(A_1)]=0$ and, therefore, the product 
$\langle S_1,\dots,S_n\rangle$ is trivial. We arrives at the contradiction, 
which proves Lemma \ref{lem2.4.5B} and Theorem  B.

{\bf Lemma \ref{lem2.4.5} C}
\label{lem2.4.5C}
{\it Under the hypothesis of Theorem C the Massey product
$\langle a\cup S_1,a\cup S_2,a\cup S_3\rangle$ is nontrivial.}

{\sl Proof.}  Let $A_1=(R(i,j))$, where $1\le i\le j\le 3$ and 
$R(i,j)\in\M_Y$,  
be a defining system for $\langle S_1,S_2,S_3\rangle$.
Then $A=(x^{j-i+1}\wedge R(i,j))$ with $1\le i\le j\le 3$ is a defining
system for $\langle a\cup S_1,a\cup S_2,a\cup S_3\rangle$. 
It is clear that $c(A)=x^3\wedge c(A_1)$, i. e.,
$[c(A)]=a^3\cup [c(A_1)]$.

The indeterminacy of a triple Massey product is simple (see 
Proposition \ref{pr1.3.9}) and we conclude that any element of
$u\in\langle a\cup S_1,a\cup S_2,a\cup S_3\rangle$  
has the form
\begin{equation}
u=[c(A)]+H\cup a\cup S_3+a\cup\overline{S}_1\cup K
\label{2.4.5}
\end{equation}
where  $H,K\in N(H^+(\widetilde{V}))$ are arbitrary matrices with appropriate
sizes and multidegrees.

Since  $H^{\ast}(\widetilde{V})$ is a free  $H^{\ast}(Y)$-module with
the basis $1, a, \dots,a^{m-1}$, expanding $a\cup H$  and $a\cup K$ in this
basis we derive from (\ref{2.4.5}) that 
\begin{equation}
u=a^3\cup[c(A_1)]+\left(\sum_{i=0}^{m-1}a^i\cup H_i\right)\cup S_3+
\overline{S}_1\cup\left(\sum_{i=0}^{m-1}a^i\cup K_i\right).
\label{2.4.6}
\end{equation}

We have  $S_i\in N(H^+(Y))$ and $m\ge 4$ and hence $u=0$ together 
with (\ref{2.4.6}) implies 
$$
[c(A_1)]+H_3\cup S_3+\overline{S}_1\cup K_3=0.
$$
This means that $0\in\langle S_1,S_2,S_3\rangle$, which contradicts to
the initial hypothesis. Hence the lemma and Theorem C are proved.

{\bf Lemma \ref{lem2.4.5} D}
\label{lem2.4.5D}
{\it Under the hypothesis of Theorem D the Massey product
$\langle a\cup S_1,a\cup S_2,a\cup S_3,a\cup S_4\rangle$ is nontrivial.}

{\sl Proof.} Suppose that the converse is true, i.e., 
$$
0\in\langle a\cup S_1,a\cup S_2,a\cup S_3,a\cup S_4\rangle.
$$
Let $A=(X(i,j))$, with $1\le i\le j\le 4$, be a defining
system for the product. 
Represent the elements of the defining system in the minimal model
$\M$ for $\widetilde{V}$ as polynomials in  $x$.
We have 
\begin{equation}
\begin{array}{c}
X(i,i)=x\wedge R_i,\ \ \ \ [R_i]=S_i,\ \ 1\le i\le 4, \\
X(i,j)=\sum_{l=0}^{K(i,j)}x^l\wedge X_l(i,j),\ \ 1\le j-i\le 2.
\end{array}
\label{2.4.7}
\end{equation}

Notice that the hypothesis of the theorem guarantees that
$\sdeg X(i,j)<2m-1<2\codim Y$, and, therefore, the generator
$y\in\M$ does not come into the expansions (\ref{2.4.7}), which means that
$X_l(i,j)\in\M_Y$ for all $i,j,l$.

For the defining system $A$ we derive the following equalities
\begin{equation}
\left\{\begin{array}{l}
1)\ \ d\,X_l(i,i+1)=\overline{R}_i\wedge R_{i+1}, \ \ l=2,\ \ i=1,2,3,  \\
2)\ \ d\,X_l(i,i+1)=0, \hskip1.8cm l\ne 2,\ \ i=1,2,3, \\
3)\ \  d\,X_{l+1}(i,i+2)=\overline{R}_i\wedge X_l(i+1,i+2)+ \\
\ \ \ \ \overline{X_l(i,i+1)}\wedge R_{i+2}, l\ge 0,\ \ i=1,2, \\
4)\ \ d\,X_0(i,i+2)=0,\ \ i=1,2.
\end{array}\right.
\label{2.4.8}
\end{equation}

This implies that $A_1(X_{j-i+1}(i,j))$, with $1\le i\le j\le 4$,
is a defining system for $\langle S_1,S_2,S_3,S_4\rangle$.

Finally we obtain      
$$
\begin{array}{c}
c(A)=\overline{X(1,1)}\wedge X(2,4)+\overline{X(1,2)}\wedge X(3,4)+
\overline{X(1,3)}\wedge X(4,4)= \\
=\sum_{l=0}^3 x^l\wedge P_l+x^4\wedge\{\overline{X_1(1,1)}\wedge X_3(2,4)+
\overline{X_2(1,2)}\wedge X_2(3,4)+ \\
+\overline{X_3(1,3)}\wedge X_1(4,4)+
\overline{X_0(1,2)}\wedge X_4(3,4)+\overline{X_1(1,2)}\wedge X_3(3,4)+ \\
+\overline{X_3(1,2)}\wedge X_1(3,4)+\overline{X_4(1,2)}\wedge X_0(3,4)\}+
\sum_{l\ge 5} x^l\wedge P_l= \\
=\sum_{l\ne 4} x^l\wedge P_l+x^4\wedge\{c(A_1)+
\overline{X_0(1,2)}\wedge X_4(3,4)+ \\
+\overline{X_1(1,2)}\wedge X_3(3,4)+\overline{X_3(1,2)}\wedge X_1(3,4)+
\overline{X_4(1,2)}\wedge X_0(3,4)\}
\end{array}
$$
where $P_l$, $l\ne 4$, are cocycles from $\M_Y$.  The cocycles $P_l$ are
expressed in terms of  $X_r(i,j)$ but for $l \ne 4$ 
their explicit forms do not matter.

Since $\sdeg c(A)<2m$, $c(A)$ is exact if and only if the cocycles $P_l$
are exact for $l=0,1,\dots$. In particular, we have
$$
\begin{array}{c}
P_4=c(A_1)+\overline{X_0(1,2)}\wedge X_4(3,4)+
\overline{X_1(1,2)}\wedge X_3(3,4)+ \\
\overline{X_3(1,2)}\wedge X_1(3,4)+
\overline{X_4(1,2)}\wedge X_0(3,4)=d\,Q,\ \ Q\in\M_Y.
\end{array}
$$
It follows from the equality 2) of the system (\ref{2.4.8}) that
no one of elements $X_l(i,i+1)$ is a cocycle for $l\ne 2$ and $i=1,2,3$.
Passing to the cohomology classes we obtain
$$
[c(A_1)]=-([\overline{X_0(1,2)}]\cup[X_4(3,4)]+
[\overline{X_1(1,2)}]\cup[X_3(3,4)]+ \\
$$
$$
+[\overline{X_3(1,2)}]\cup[X_1(3,4)]+[\overline{X_4(1,2)}]\cup[X_0(3,4)]).
$$
Therefore, $[c(A_1)]\in N(H^+(Y)\cdot H^+(Y))$ and that contradicts to the 
irreducibility of   
$\langle S_1,S_2,S_3,S_4\rangle$. This proves the lemma and Theorem D.

{\sl Remark.}  The matrix cohomology classes  
$[X_l(1,2)]$ and $[X_l(3,4)]$ where $l=0,1,3,4$ are not arbitrary.
The equalities (\ref{2.4.8}) imply relations   
$$
\begin{array}{cc}
\overline{S}_1\cup [X_l(2,3)]+[\overline{X_l(1,2)}]\cup S_3=0, & \\
                                                               & l=0,1,3,4. \\
\overline{S}_2\cup [X_l(3,4)]+[\overline{X_l(2,3)}]\cup S_4=0, & \\
\end{array}
$$
The classes $[X_l(2,3)]$ with $l=0,1,3,4$ are restricted only by
the degrees of their matrix entries.

\end{document}